\documentclass[a4paper,11pt]{article}

\def\thebibliography#1{\section*{Literature Cited\markboth
 {Literature Cited}{Literature Cited}}\list
 {[\arabic{enumi}]}{\settowidth\labelwidth{[#1]}\leftmargin\labelwidth
 \advance\leftmargin\labelsep
 \usecounter{enumi}}
 \def\newblock{\hskip .11em plus .33em minus -.07em}
 \sloppy
 \sfcode`\.=1000\relax}

\usepackage{amssymb}
\usepackage{mathrsfs}
\usepackage{graphicx}
\usepackage{amsmath}
\usepackage{amsfonts}
\usepackage{amsthm}
\usepackage{epic}
\usepackage{eepic}

\newtheorem{Th}{Theorem}
\newtheorem{Propn}{Proposition}

\newtheorem{Lmm}{Lemma}

\begin{document}
\begin{center}

 {\LARGE 
On the Inviscid Burgers Equation and the Axiom of Choice}\\ 		
 \setlength{\baselineskip}{3ex}

\vspace{.1in}
		   
John M. Noble\\ 
Matematiska institutionen,\\
Linköpings universitet,\\
58183 LINK\"OPING, Sweden\\

  \end{center}

\paragraph{Summary} This article gives an example where  
the use of Tychonov's theorem, which is equivalent to the Axiom of Choice,
yields two contradictory results. That is, by using Tychonov's
theorem, full proofs of contradictory results are obtained. That the
Choice Axiom implies Tychonov's theorem is standard. That Tychonov's
theorem implies the Axiom of Choice is a result of Kelley, found in the
article~\cite{Ke}. 

This paper discusses the first direction; Choice implies Tychonov's
theorem. The example discussed in this article provides a counter example to the
assertion that the unit ball in
$L^2$ is   relatively weakly compact.     
Tychonov's theorem states that  the unit ball in
$L^2$ is relatively weakly compact. The proof of Tychonov's theorem crucially uses Choice and Kelley shows that Choice is equivalent to Tychonov's theorem. Therefore, as a corollary of the proof that contradictory results may be obtained by assuming the relative weak compactness of the unit ball in $L^2$, it follows that the Choice Axiom is inadmissible  in
mathematical analysis.

The article is structured as follows. Section (\ref{intro})  is 
introductory. Section (\ref{downward}) illustrates the
construction of solutions to the inviscid Burgers equation in
terms of the velocities of Euler Lagrange trajectories. It also computes a
formula for the evolution of downward jumps for solutions to the inviscid
Burgers' equation. 
Section (\ref{varadhansect})  sketches the large deviations argument (for  smooth, bounded potentials) which leads to the
representation of the solution to the inviscid Burgers equation in terms
of the velocity of the trajectory that {\em minimises} the associated action
functional, subject to the appropriate constraints. Also outlined are the standard arguments from
the Calculus of Variations that prove the existence of a trajectory at which
the minimum of the action functional is attained subject to the appropriate
constraints, that this trajectory solves the associated Euler Lagrange equations, and that the  solutions
of the inviscid Burgers' equation may be constructed using these  
minimising trajectories. 

Section
(\ref{final})  presents an example of a smooth, bounded, space / time
periodic potential, for which the viscosity solution to the inviscid Burgers'
equation with that potential may {\em not} be constructed using  
trajectories that {\em  minimise} the associated action functional.  The
solution may be constructed using {\em critical points} of the associated
action functional, but these critical points are {\em not} the
minimising trajectories. \clearpage
 
\section {Introductory}~\label{intro}  

Let $W = C_0({\bf R}_+)$; that is, continuous functions $f: {\bf R}_+ \rightarrow {\bf R}$ such that $f(0) = 0$. Let 
$w$ denote a trajectory for standard Brownian motion in
${\bf R}$ satisfying $w(0) = 0$. Let \[ (W, {\cal F}, ({\cal F}_{s,t})_{0 \leq s
\leq t < +\infty}, {\bf P}) \] denote the filtered Gaussian probability space associated
with
$w$, where
${\cal F}_{s,t}$ is the Borel sigma algebra generated by the increments of continuous functions between $s$ and $t$; $(w(u) -
w(v))_{s
\leq v
\leq u \leq t}$, ${\cal F} = \cup_{0 \leq s \leq t < +\infty} {\cal F}_{s,t}$, and ${\bf P}$ is the probability measure associated with standard Brownian motion. That is, under ${\bf P}$, $w(0) = 0$ with probability $1$, for any collection $(t_1, \ldots, t_{n+1})$ with $0 \leq t_1 \leq \ldots \leq t_{n+1} < +\infty$ the random variables $(w(t_{j+1}) - w(t_j))_{j=1}^n$ are independent Gaussian random variables with $w(t_{j+1}) - w(t_j) \sim N(0, t_{j+1} - t_j)$ and, for all $ s < t$, $s \leq u \leq v \leq t$, $w(v) - w(u)$ is ${\cal F}_{s,t}$ measurable. 

Let
$E_{\bf P}$ denote the expectation operator with respect to ${\bf P}$. Let
$V$ denote a smooth, bounded (time dependent) potential $V : {\bf R}_+ \times
{\bf R}
\rightarrow  {\bf R}$.  Using subscripts to denote derivatives with
respect to the subscripted variable, consider the equation

\[ \left\{\begin{array}{l} U^{(\epsilon)}_t =
\frac{\epsilon}{2}U^{(\epsilon)}_{xx} -
\frac{1}{\epsilon}U^{(\epsilon)}V \\ U^{(\epsilon)}(0,.) =
e^{-\frac{1}{\epsilon}\phi(.)}. \end{array}
\right.\]

\noindent A Feynman - Kacs representation of the solution may be
employed;

\begin{equation}\label{fkexp} U^{(\epsilon)}(t,x) = E_{\bf P} \left [
\exp\left\{-\frac{1}{\epsilon}\left(\phi(x+\sqrt{\epsilon}w(t)) + \int_0^t
V(t-s, x + \sqrt{\epsilon}w(s))ds
\right)\right\}\right].\end{equation}

\noindent Consider now $v^{(\epsilon)} = -\epsilon \log U^{(\epsilon)}$ and
note that $v^{(\epsilon)}$ satisfies 

\[ \left\{\begin{array}{l} v_t^{(\epsilon)} =
\frac{\epsilon}{2}v_{xx}^{(\epsilon)} -
\frac{1}{2}(v^{(\epsilon)}_x)^2 + V \\ v^{(\epsilon)}(0,.) = \phi(.).
\end{array}\right.\]

\noindent From the Feynman Kacs representation for $U^{(\epsilon)}$,
given in equation (\ref{fkexp}),  
$v$ satisfies

\[ v^{(\epsilon)}(t,x) = - \epsilon \log E_{\bf P} 
\left [\exp \left \{ -
\frac{1}{\epsilon}\left(\phi(x + \sqrt{\epsilon} w(t)) +  \int_0^t V(s,x +
\sqrt{\epsilon}w(t-s)) ds\right) \right
\}
\right ].\]

\noindent  Suppose that $\phi$ and $V$ are uniformly bounded. Let ${\cal
S}_n$ denote the space of functions $f : {\bf R}\rightarrow {\bf R}$, with
bounded derivative $\dot{f}$, such that $\dot{f}$ is piecewise constant on
intervals
$[\frac{k}{2^n},
\frac{k+1}{2^n})$; that is, if
$f
\in {\cal S}_n$ then there exist real numbers $(\lambda_j)_{j=-\infty}^\infty$
such that 

\[ \dot{f} = \sum_j \lambda_j \chi_{[\frac{j}{2^n},
\frac{j+1}{2^n})},\]

\noindent where $\chi_A$ denotes the indicator function for a set $A$. This notation will be used throughout. Let ${\cal S} = \cup_n {\cal S}_n$. Then, by  Varadhan's theorem  from large deviations, it follows that 

\begin{equation}\label{varadhan}
v(t,x) := \lim_{\epsilon \rightarrow 0} v^{(\epsilon)}(t,x) = \inf_{\xi :
\xi(t) = x,\; \xi \in {\cal S}}\left \{\phi(\xi(0)) + 
\frac{1}{2}\int_0^t \dot{\xi}^2(s) ds + \int_0^t V(s,\xi(s)) ds\right \} .
\end{equation}

\noindent  The result is well known, but a full proof   is given in section (\ref{varadhansect}), for completeness of the presentation. The treatment is taken from Dembo and 
Zeitouni~\cite{DZ}, with appropriate simplifications, because only a special case of their setting is required here. This approach, of taking a partition over the time interval and letting the mesh size tend to zero, seems more appropriate in a situation where it is convenient to locate where Choice enters. Choice is used in the proof, but it seems clear that the statement of the result can be modified and the proof modified to produce a   version of the result where the proof does not require Choice. Another full and elegant treatment  
may be found in Deuschel and Stroock~\cite{DS}.   

The next question is whether or not there exists a
minimiser in the space $W^{1,2}({\bf R})$; namely, a trajectory at which the
infimum is attained. For this specific problem, there exists a trajectory at which the minimum of the action functional is attained. This is a 
consequence of the relative weak compactness of a ball of finite radius in
$L^2$. The result was first established by   Tonelli in~\cite{T1} and~\cite{T2}. It is treated in a more
general framework by 
 Cesari in~\cite{Ces}. The relative weak compactness of the unit ball in $L^2$ is crucial here. 

The whole point at issue seems to be the existence of the minimising trajectory, because under the conditions in this article, having established existence of  minimising trajectories, it is relatively straightforward to show that they satisfy the associated Euler Lagrange equations. 

 In more
general problems within the Calculus of Variations, where the action functional does not satisfy
some crucial hypotheses of the framework of Tonelli,  there may not exist minimising trajectories, as may be seen, for example, in
Ball and Mizel~\cite{B-M} and~\cite{BHJPS}, where examples are given of
situations where global minimisers exist, but do not  satisfy the Euler
Lagrange equations and examples of situations where global minimisers do
not exist. 

In the
example presented in section (\ref{final}) of this article,   the resulting process  
does indeed satisfy the Euler Lagrange equations in the limit, but it is shown that the
trajectory picked out is not   the minimiser. Therefore,   the minimiser is
not necessarily the trajectory that appears in representation of the
`viscosity' solution of the inviscid Burger equation using solutions to the
associated Euler Lagrange equations. \vspace{5mm}

\noindent {\bf Observation:} It  is the {\em countable}
version of the Choice Axiom that is shown to lead to contradictory results in this
article. \vspace{5mm}

\noindent Let $u^{(\epsilon)} = v^{(\epsilon)}_x$, then $u^{(\epsilon)}$
satisfies

\begin{equation}\label{burvis} \left\{\begin{array}{l}u^{(\epsilon)}_t +
\frac{1}{2}(u^{(\epsilon)2})_x =
\frac{\epsilon}{2}u^{(\epsilon)}_{xx} + V_x\\
u^{(\epsilon)}(0,.) = \phi_x(.).\end{array}\right.\end{equation}

\noindent Let $u = \lim_{\epsilon \rightarrow 0} u^{(\epsilon)}$, where
the limit is taken in the relative weak topology in $L^2$. Then the limit
exists and satisfies 

\begin{equation}\label{invbur} \left\{\begin{array}{l}u_t +
\frac{1}{2}(u^{2})_x =
 V_x\\
u (0,.) = \phi_x(.).\end{array}\right.\end{equation}

\noindent Let 

\begin{equation}\label{action} {\cal A}(\xi;t,x) = \frac{1}{2}\int_0^t
\dot{\xi}^2(s)ds +
\int_0^t V(s,
\xi(s))ds + \phi(\xi(0)).\end{equation}

\noindent If $V$ is infinitely differentiable and uniformly bounded (as it is in the example considered in this article), then it is standard that for each $(t,x)$ there
exists a trajectory
$\xi^{(t,x)}$ that minimises ${\cal A}$ subject to the constraint that
$\xi^{(t,x)}(t) = x$. Throughout, the following notation will be used: for a function $f$ of two arguments, $\dot{f}$ will be used to denote the derivative with respect to the first argument and $f^\prime$ will be used to denote the derivative with respect to the second argument. If $f$ has only one argument, then either $f^\prime$ or $\dot{f}$ may be used to denote the derivative with respect to that argument. Usually $\dot{f}$ will be employed if the argument serves as a `time' variable and $f^\prime$ will be used if the argument serves as a space variable. It follows that $V^\prime = V_x$ and $\phi^\prime = \phi_x$. For the action functional defined by equation (\ref{action}), relatively straightforward arguments using the relative weak compactness of the unit ball, show that there exists a trajectory where the global minimum of the action functional is attained and, having established existence, relatively straightforward arguments from the calculus of variations show that any minimiser satisfies the Euler Lagrange equations, 

\begin{equation}\label{elnew2}\left\{\begin{array}{l} \ddot{\xi}^{(t,x)}(s) =
V^\prime (s,\xi^{(t,x)}(s))
\\
\xi^{(t,x)}(t) = x, \qquad \dot{\xi}^{(t,x)}(0) = \phi^\prime ( \xi^{(t,x)}(0)).
\end{array}\right. \end{equation}

\noindent In general, solutions to equation (\ref{elnew2}) will not
be unique if $t$ is sufficiently large. Having established
existence of trajectories where the global minimum is attained and that they solve equation
(\ref{elnew2}), it may be shown that the solution to equation (\ref{invbur})
has representation in terms of the minimising trajectories. Let $\eta$
solve

\begin{equation}\label{elnew}
\left\{\begin{array}{l} \ddot{\eta}  = V^\prime (t,\eta ) \\
\eta(0,x) = x, \qquad \dot{\eta}(0,x) = \phi_x (x),\end{array}\right.
\end{equation}

\noindent then $u$, the solution to equation (\ref{invbur}) has a 
representation

\begin{equation}\label{rep} u(t,x) = \dot{\eta} (t,
\eta^{-1}(t,x)) = \dot{\xi}^{(t,x)}(t),\end{equation}

\noindent where $\xi$ solves equation (\ref{elnew2}). It turns out that $\xi$ also minimises the action functional (\ref{action}) subject to the constraint that $\xi(t) = x$. These results are all standard and their proofs are all outlined in the article.
\vspace{5mm}

\noindent The counter example in section (\ref{final}) gives an example
where  the viscosity solution to equation (\ref{invbur}) may be constructed
using trajectories which satisfy the Euler Lagrange equations, but that
these trajectories {\em necessarily} do not minimise the action functional. The
strategy is as follows: The potential chosen is

\[ V(t,x) = \cos (\sin(t)) - \cos (x + \sin(t)).\]

\noindent The time dependence seems to be crucial here to manufacture a counter example. For fixed initial condition, there is uniqueness of solution to equation (\ref{burvis}). With this choice of potential, it is shown that there exists
exactly one initial condition $\phi^{(\epsilon)}_x$ that yields a {\em periodic} solution to equation (\ref{burvis}).  It is shown that these periodic solutions $u^{(\epsilon)}$ converge (in the relative weak topology) to a limit $u$ and that this is the only {\em periodic} solution of equation
(\ref{invbur}); $\phi_x$, the weak limit of $\phi^{(\epsilon)}_x$ is the only initial condition that yields periodic solutions to equation (\ref{invbur}). It is shown that there is
necessarily a {\em periodic modulo $2\pi$} solution to the Euler Lagrange equations
involved in the construction of this periodic solution to equation
(\ref{invbur}) and, indeed, that all trajectories used in the construction,
when run backwards, converge to a {\em periodic} trajectory. But it is
shown that the only two {\em periodic modulo $2\pi$} trajectories  which
solve equation

\[ \ddot{\xi} = \sin (\xi + \sin (t))\]

\noindent are 

\[ \xi (t) = t - \sin(t)\]

\noindent and

\[ \xi (t) = \pi - t - \sin(t).\]

\noindent {\bf Note} Periodic modulo $2\pi$ means that $\xi(t) \mod(2\pi)$ is periodic. It is shown that, for sufficiently large $t$, neither of these
minimise the action functional

\[ {\cal A}(\xi; t) := \frac{1}{2}\int_0^t \dot{\xi}^2(s)ds + \int_0^t \cos
(\sin(s))ds - \int_0^t \cos (\xi(s) + \sin(s))ds + \phi( \xi(0)).\]

\noindent The problem became apparent 
following results in  two articles by E, Khanin, Mazel and
Sinai~\cite{EKMS1} and~\cite{EKMS2}. The article~\cite{N} presents 
analysis of the moments of the stochastic inviscid Burgers' equation, under
a special case of the hypotheses considered in   E, Khanin, Mazel, Sinai, 
which are of interest following the invariant measure proved in
~\cite{EKMS2}. The article~\cite{EKMS2} used  crucially the existence of  
 trajectories where   the minimum of the action functional is attained, that  minimising trajectories satisfy the associated
Euler Lagrange equations (\ref{elnew}), and that it is minimising trajectories that are  used to construct the solution given in equation (\ref{rep}) to equation
(\ref{invbur}). Conditioned on this step, every other part of the argument in
the article~\cite{EKMS2} is clear.

\section{The Inviscid Burgers Equation and the Euler Lagrange
Equations} \label{downward}  
In general, for the inviscid Burgers' equation with smooth initial condition and smooth potential, there
will be uniqueness of solutions to the associated  Euler Lagrange equations with appropriate
boundary data up to the onset of downward jumps. Up to the onset of
downward jumps, there will be exactly one relevant Euler Lagrange
trajectory and this will be the global minimiser. After the onset of
downward jumps, there will be a family of solutions to the Euler
Lagrange equation with relevant boundary data. The downward jumps must
evolve in such a way that the inviscid Burgers equation is satisfied.

Since the material in section
(\ref{downward}) is standard, the proofs are only sketched. 

\paragraph{The Representation in terms of Euler Lagrange Trajectories}
Consider the equation

\begin{equation}\label{youuu} \left\{\begin{array}{l} u_t +
\frac{1}{2}(u^2)_x = V_x
\\ u_0 =
\phi_x
\end{array}\right.\end{equation}

\noindent where $\phi$ is periodic and Lipschitz, and $\phi_{xx}$ is
bounded from above, and
$V$ is smooth and periodic, with all derivatives uniformly bounded. Let
$\theta$ satisfy

\begin{equation}\label{theteqn} \left\{\begin{array}{l} \dot{\theta} =
u(t,\theta )
\\
\theta (0,x) = x,
\qquad
\dot{\theta}(0,x) = \phi_x(x).
\end{array}\right.\end{equation}

\noindent Let $u^\prime $ denote the function such that $u^\prime (t,x) =
u_x(t,x)$ and let $\dot{u}$ denote the function such that $\dot{u}(t,x) =
u_t(t,x)$ and $V^\prime$ the function such that $V^\prime (t,x) = V_x(t,x)$,
$\phi^{\prime \prime}$ such that $\phi^{\prime \prime}(x) = \phi_{xx}(x)$
and $V^{\prime \prime}$ such that $V^{\prime \prime}(t,x) = V_{xx}(t,x)$. In short, for a function $f \in C^\infty ({\bf R}\times {\bf S}^1)$, where ${\bf S}^1$ is used to denote the circle $[0, 2\pi)$ where $0$ is identified with $2\pi$, $\dot{f}$ denotes the derivative with respect to the first argument and $f^\prime$ denotes the derivative with respect to the second argument. If $f$ is a function of only one variable, then either $\dot{f}$ or  $f^\prime$ will be used to denote the derivative.

\noindent Set $u(t,x-) := \lim_{y \uparrow x}u(t,y)$ and
$ u(t,x+) = \lim_{y \downarrow x}u(t,y)$. Then it follows that for all $x$,
$u(t,x-) \geq u(t,x+)$ and, in particular, $u(t,\theta(t,x)-) \geq
u(t,\theta(t,x)+)$. \vspace{5mm}

\noindent Let $\eta$ satisfy 

\begin{equation}\label{etaeqn}\left\{\begin{array}{l}\ddot{\eta} = V^\prime
(t,\eta)\\
\eta(0,x) = x,
\qquad
\dot{\eta}(0,x) = \phi_x(x).\end{array}
\right.\end{equation}

\begin{Lmm} \label{dervub} Let  $u$ satisfy equation (\ref{youuu}). Suppose that $\sup_{s,x}|V^\prime (s,x)| + \sup_{s,x} |V^{\prime \prime}(s,x)| < +\infty$ and $\sup_x \phi^{\prime \prime}(x) < +\infty$ and $\sup_x |\phi^\prime(x)| < +\infty$. Then
\[ \sup_{0 \leq s \leq t} \sup_x u_x(s,x) < C(t) < +\infty,\]
\noindent where
 \begin{equation}\label{setcee}
 C(t) = \sup_x  \phi^{\prime \prime} (x)  + t \sup_{s,x}
 V^{\prime \prime}(s,x).
\end{equation}
 That is, there is an {\em upper bound} on the derivative. 
\end{Lmm}

\noindent {\bf Proof} Set $w = u_x$. Then $w$ satisfies

\[\left\{\begin{array}{l} w_t + w^2 + uw_x = V_{xx} \\
w(0,x) = \phi_{xx}.\end{array}\right.\]

\noindent Let $\tilde{\theta}$ satisfy

\begin{equation}\label{tildetheta}\left\{ \begin{array}{ll} \dot{\tilde{\theta}}^{(t)}(s,x) = - u(t-s, \tilde{\theta}^{(t)}(s,x)) & s \geq 0 \\ \tilde{\theta}^{(t)}(0,x) = x & 
\end{array}
 \right.\end{equation}

\noindent It is clear that 

\[ w(t,x) = \phi^{\prime \prime}(\tilde{\theta}^{(t)}(t,x)) + \int_0^t V^{\prime
\prime} (s, \tilde{\theta}^{(t)}(t-s,x))ds - \int_0^t w^2 (s, \tilde{\theta}^{(t)}(t-s,x))ds.
\]

\noindent From   the hypotheses on
$\phi$ and
$V$, it follows that  for any $t < +\infty$, 

\[ \sup_{0 \leq s \leq t}\sup_{x \in [0,2\pi)}w(s,x) \leq C(t) < +\infty\]

\noindent and the result follows. \qed \vspace{5mm}

\noindent Since the derivative is bounded from above, there can be no `upward jumps' in a solution; any discontinuities have to be `downward jumps'. The function of the next lemma is to show that when Euler Lagrange trajectories are being used to construct the solution, a trajectory is used until it enters a downward jump. After this, it is no longer used in the construction. \vspace{5mm}

\begin{Lmm}\label{disappear} Recall that ${\cal S}^1$ denotes the circle $[0,2\pi)$, with the identification $0 = 2\pi$. Set 

\[ {\cal S}(t) = \{x \in {\bf S}^1 | \theta_x(t,x) = 0\}.\]

\noindent  Then, for all $(s,t)$ such that $s \leq t$, ${\cal S}(s) \subseteq {\cal S}(t)$. 
\end{Lmm}

\noindent {\bf Proof} Recall the notation $w = u_x$ and set 

\[ f(t,x) =    
w(t,\theta(t,x)) .   \] 

\noindent Recall that
$\dot{\theta} = u(t,\theta)$. Since $\theta(0,x) = x$, it follows that $\theta_x(0,x) \equiv 1$, so that, directly from equation (\ref{theteqn}), together with the definition of $f$, 

\begin{equation}\label{dottheta}\left\{\begin{array}{l} \dot{\theta}_x = f(t,x)\theta_x, \\\theta_x(0,x) \equiv 1,\end{array}\right. 
\end{equation}

\noindent for all $t \geq 0$ and all $x \in {\bf S}^1$. Set $\sigma(x) = \inf\{t | \theta_x(t,x) = 0\}$. It now follows directly from the upper bound in lemma (\ref{dervub})  that $\theta_x(t,x) \equiv 0$ for all $t > \sigma(x)$. Furthermore, for all $0 \leq s \leq t < +\infty$, 

\begin{equation}\label{thetaxbd} 0 \leq \theta_x(t,x) = \theta_x(s,x) \exp\left\{ \int_s^t f(r,x)dr \right\} < \theta_x(s,x) \exp\{ (t-s) C(t)\} < +\infty, 
\end{equation}

\noindent from which the result follows directly.
\qed
\vspace{5mm}

\noindent  For all $t > \sigma(x)$,  
$\theta_x (t,x) = 0$. Recall that ${\cal S}(t)= \{x | \theta_x(t,x) = 0\}$ and set ${\cal D}(t) = {\bf S}^1 \backslash {\cal S}(t)$. For $y \in \theta (t,{\cal D}(t))$, note that 
\[ u(t,y) = \dot{\theta}(t, \theta^{-1}(t,y))\]

\noindent and note that $|\theta(t, {\cal D}(t))|:= \int_{{\cal D}(t)} \theta_x(t,x) dx = 2\pi$. The sets ${\cal D}(t)$ and $\theta (t,{\cal D}(t))$ are open. For $0 \leq t \leq \sigma(x)$, $\theta$ satisfies

\[ \ddot{\theta}(t,x) = \frac{d}{dt} u(t,\theta(t,x)) = \dot{u}(t,\theta(t,x)) + \dot{\theta}(t,x) u^\prime (t,\theta(t,x)) = V(t,\theta(t,x)).\]

\noindent It follows that $\theta(t,x) = \eta(t,x)$ for all $0 \leq s \leq t$ where $\eta$ satisfies equation  (\ref{etaeqn}). It follows that, for $x \in \theta (t,{\cal D}(t))$,

\begin{equation}\label{construct} 
 u(t,x) = \dot{\eta}(t,\eta^{-1}(t,x)),
\end{equation}

\noindent where $\eta$ satisfies equation (\ref{etaeqn}). Furthermore, for $x \in \theta (t,{\cal D}(t))$, the inverse $\eta^{-1}(t,x)$ is uniquely defined. A càdlàg version of the solution is given by equation (\ref{construct}) for $x \in \theta (t,{\cal D}(t))$ and $u(t,x) = \lim_{y \downarrow x, y \in  \theta (t,{\cal D}(t))} u(t,y)$ for $x \in [0,2\pi)$. Therefore, a trajectory 
$\eta(t,x)$ which solves equation (\ref{etaeqn}) is used in the construction
for all 
$t \in (0,\sigma(x))$, but is not used for any $t > \sigma(x)$. \vspace{5mm}

\noindent The following computation shows how the downward jumps evolve. After the `onset of downward jumps'; namely, for $t > T$
where 
$T = \inf_x \sigma(x)$, the set ${\cal S}(x) = \{y |
\eta(t,y) = x\}$ may contain more than one element. But the choice of
trajectories that may be used in the construction of the solution is not
arbitrary. The equation (\ref{youuu}) determines how the downward jump
sites must evolve.

\paragraph{The Onset of Downward Jumps} Consider the Inviscid Burgers'
Equation (equation (\ref{invbur})). Suppose that $\phi$ is $2\pi$
periodic and
$V$ is $2\pi$ periodic in both variables, with $\int_0^{2\pi}\phi(x)dx
= 0$ and
$\int_0^{2\pi}V(t,x)dx = 0$ for each $t \geq 0$. Then it is straightforward to compute that

\[ u(t,x) = \phi_x(\tilde{\theta}^{(t)}(t,x)) + \int_0^t V^\prime (s, \tilde{\theta}^{(t)}(t-s, x))ds,\]

\noindent where $\tilde{\theta}$ satisfies equation (\ref{tildetheta}). It follows that 

\[ \sup_{0 \leq s \leq t} \sup_{0 \leq x \leq 2\pi} |u(s,x)| \leq \sup_x
|\phi_x(x)| + t \sup_{0 \leq s \leq t}\sup_{0 \leq x \leq 2\pi}|V_x(s,x)|.\]

\noindent Suppose that $u(t,.)$ has a `downward
jump' at site
$\theta (t)$. The following analysis shows how the position of the downward
jump  evolves in time.  Since  $u \in L^\infty $, consider integration
against  test functions
$\psi \in L^1$. 

\begin{eqnarray}\nonumber \int_0^\infty \int_{{\bf R}}
\psi(s,y)u_s(s,y)dyds &+& 
\int_0^\infty \int_{\bf R}
\frac{1}{2}\psi(s,y) (u^2)_x(s,y) dy ds\\  \label{weak} &=& \int_0^\infty 
\int_{\bf R}
\psi(s,y) V_x(s,y)dyds.\end{eqnarray}\vspace{5mm}

\begin{Lmm}[Downward Jump Evolution]\label{down2}If a downward jump
develops,  found at
site $\theta(t)$ at time $t$, then $\theta(t)$   evolves according to the
equation

\begin{equation}\label{down}
 \dot{\theta}(t) =  \frac{u (t,\theta(t)+) + u (t,\theta(t)-)}{2}.
\end{equation}
\end{Lmm}\vspace{5mm}

\noindent {\bf Proof (sketch) of lemma (\ref{down2})} The calculation leading to
 formula (\ref{down}) is now outlined.  Assume that a downward
jump in the space variable develops at site $\theta(t)$, time $t$ and
continues along trajectory $\theta(s)$ for $s \geq t$.    Consider  test
functions
$\psi^{(\delta)} \in C^\infty$, such that $\psi^{(\delta)}(s,x) = 0$
for all $s < t$, $\sup_{0 <
\delta < 1}\sup_{(s,x) 
\in {\bf R}_+\times {\bf R}}|\psi^{(\delta)}(s,x)| < K < +\infty$ for some
$K$ and   with support
within a tube of radius $\delta$ around the graph
$(s,\theta(s))_{s \geq t}$. Suppose, furthermore, that $\psi^{(\delta)}$ are chosen in such a way that there exists a
function
$f\in C^\infty ({\bf R})$ such that $f(s) = 0$ for $s \leq t$ and 
$\psi^{(\delta)}(s,\theta(s)) = f(s)$ for all
$\delta > 0$ and
$s
\geq t$.  Then, from equation (\ref{weak}) it follows that

\begin{eqnarray*}\lefteqn{ \int f (s(\theta)) (u(s(\theta), \theta-) -
u(s(\theta), \theta+)) d\theta}
\\&& \hspace{20mm}  + \frac{1}{2}\int
 f (s )(u^2(s,\theta(s)+) - u^2(s,\theta(s)-))ds = 0 
\end{eqnarray*}

\noindent so that, for any test function $f \in C^\infty ({\bf R})$,  

\begin{eqnarray*}\lefteqn{
\int f(s) (u(s , \theta(s)-) - u(s , \theta(s)+)) \dot{\theta}(s) ds}
\\&& \hspace{20mm}  + \frac{1}{2}\int
f(s)(u^2(s,\theta (s)+) - u^2(s, \theta (s)-))ds = 0, 
\end{eqnarray*}

\noindent   Equation
(\ref{down}) follows, thus proving lemma (\ref{down2}).\qed
\vspace{5mm}

\section{The Large Deviations and Calculus of Variations
Arguments}\label{varadhansect}  

Suppose 
$V \in C^\infty ({\bf R}_+ \times {\bf R})$ is smooth, bounded and periodic
in both variables. This section  outlines the  standard and well known
arguments to show that the viscosity solution to the inviscid Burgers
equation may be represented in terms of the velocity of solutions to the
associated Euler Lagrange equations, that these solutions are the  critical
points of the associated action functional and that they are minimising trajectories 
of the action functional.  \vspace{5mm}

\noindent Let
$U^{(\epsilon)}(t,x; \phi)$ denote the solution to the equation

\begin{equation}\label{kac} \left\{\begin{array}{l}  
U^{(\epsilon)}_t = \frac{\epsilon}{2} 
U^{(\epsilon)}_{xx} -
\frac{1}{\epsilon}U^{(\epsilon)}V \\ U^{(\epsilon)}(0,x; \phi) =
\exp\left\{-\frac{1}{\epsilon} \phi (x)\right\},\end{array}\right. 
\end{equation} 

\noindent where $\phi $ is a bounded Lipschitz function.
Set  $u^{(\epsilon)} = -\epsilon\frac{\partial}{\partial x}\log
U^{(\epsilon)}$, then $u^{(\epsilon)}$ satisfies

 \begin{equation} \label{viscous}
 \left\{\begin{array}{l} u^{(\epsilon)}_t +
\frac{1}{2}(u^{(\epsilon)2})_x =
\frac{\epsilon}{2}  u^{(\epsilon)}_{xx} +  V_x \\ u(0,.) =
 \phi_x .\end{array}\right.\end{equation}

\noindent A weak limit $u$ in $L^2$ of $(u^{(\epsilon)})_{\epsilon > 0}$ will provide   a `viscosity' solution to the equation

\begin{equation}\label{burrr} \left\{\begin{array}{l} u_t +
\frac{1}{2}(u^2)_x= V_x \\
u(0,.) =  \phi_x.\end{array}\right.\end{equation}

\noindent Writing out the Feynman - Kac representation of the solution to
equation (\ref{kac}), using $E_{\bf P}$ as expectation over standard Brownian
motion with initial condition $w_0 = 0$ gives

\[ U^{(\epsilon)}(t,x) = E_{\bf P}\left[ \exp\left \{- \frac{1}{\epsilon}\left
(\int_0^t V(t-s, \sqrt{\epsilon}w_s +x) ds +
\phi (\sqrt{\epsilon}w_t+x)\right )
\right \}\right ].\]\vspace{5mm}

\noindent The following presents a special case of Varadhan's
theorem.\vspace{5mm}

\begin{Th}[Varadhan] \label{varad}Suppose that $V$ is smooth and bounded
and that
$\phi$ is Lipschitz and bounded. Let ${\cal S}_n$ denote the space of
functions $\xi \in C^1 ({\bf R})$ such that the derivative $\dot{\xi}$ satisfies  

\[\dot{\xi} = \sum_{k=1}^{2^n+1} \lambda_{n,k} \chi_{  [\frac{(k-1)t}{2^n}, \frac{kt}{2^n}  )},\]

\noindent where $(\lambda_{n,k})_{k=1 }^{2^n+1}$ is a collection of real numbers. Let ${\cal S} =
\cup_n {\cal S}_n$. Then

\[ \lim_{\epsilon \rightarrow 0} -\epsilon \log U^{(\epsilon)}(t,x) =
\inf_{\xi \in {\cal S} : \xi(t) = x}\left\{\frac{1}{2}\int_0^t \dot{\xi}^2(s) ds
+
\int_0^t V(s,
\xi(s))ds + \phi(\xi (0))
\right\}.\]
\end{Th}\vspace{5mm}

\noindent {\bf Proof} The proof of this basically follows Dembo and
Zeitouni~\cite{DZ}, with appropriate simplifications because only a
special case of their results is required here. The proof is carried out in
steps. Firstly,  set
$t_k^{(n)} =\frac{k}{2^n}t$ and let
$(Z_j)_{j=1}^{2^n}$ denote independent random variables, each with
distribution $N(0,
\frac{\epsilon t}{2^n})$ (normal, expected value $0$ and variance $\frac{\epsilon t}{2^n}$).  Set $k(s) = \left [ \frac{2^n s}{t}\right ]$, where
$[.]$ denotes the integer part, so that  

\[ k(s) = \sup \left \{k \in {\bf Z} : \frac{kt}{2^n} < s\right \}.\]

\noindent  For $z \in {\bf R}^{2^n}$, set

\begin{equation}\label{wayyen} Y^{(n)}(s,z) = x + \sum_{j=1}^{k(s)} z_j +
\left(
\frac{2^n s}{t} - k(s)\right) z_{k(s) + 1} \qquad 0 \leq s \leq t.
\end{equation}

\noindent Let ${\bf Q}^{(\epsilon,n)}$
denote the probability measure with respect to $Z = (Z_j)_{j=1}^{2^n}$.
Then, using 
$E^{(n,\epsilon)}$ to denote  expectation with respect to the probability
measure ${\bf Q}^{(\epsilon,n)}$, set   

\begin{equation}\label{youeeen}U^{(\epsilon, n)}(t,x) = E^{(\epsilon,n)}\left[
\exp\left\{-\frac{1}{\epsilon}\left(\phi(Y^{(n)}(t,Z)) +  \int_0^t V(t - s,
Y^{(n)}(s,Z))ds\right)\right\}\right ].
\end{equation}

\noindent For $\lambda \in {\bf R}^{2^n}$, set

\[ \Lambda^{(n)}(\lambda) = \epsilon \log
E^{(\epsilon,n)}\left [e^{\frac{1}{\epsilon}\sum_{j=1}^{2^n} \lambda_j Z_j} \right ] =
\frac{1}{2}
\sum_{j=1}^{2^n}
\frac{t}{2^n} \lambda_j^2. \]

\noindent  Let $\Lambda^{(n)*}$ denote the
Fenchel-Legendre transform of $\Lambda^{(n)}$; namely, for $x \in {\bf
R}^{2^n}$, set 

\begin{equation}\label{fenleg1} \Lambda^{(n)*}(x) = \sup_{\lambda \in {\bf
R}^{2^n}}
\left\{\sum_{j=1}^{2^n}
\lambda_j x_j - \Lambda^{(n)}(\lambda)\right\}\end{equation}

\noindent and note that 

\begin{equation}\label{fenlegtran} \Lambda^{(n)*}(x) =
 \frac{2^n}{2t}\sum_{j=1}^{2^n}
x_j^2. 
\end{equation}

\noindent The following result is a
simplified version of the Gärtner - Ellis theorem, following
the presentation in~\cite{DZ}. 

\begin{Th}[Gärtner - Ellis] \label{gellis} For any closed set $F \subset {\bf
R}^{2^n}$, 
\begin{equation}\label{uppbdd} \limsup_{ \epsilon \rightarrow 0}\epsilon
\log {\bf Q}^{(\epsilon,n)}\{  F\}
\leq -\inf_{x \in F} \Lambda^{(n)*}(x)
 \end{equation}

\noindent For any open set $G \subset {\bf R}^{2^n}$
\begin{equation}\label{lowbdd} \liminf_{ \epsilon \rightarrow 0}\epsilon
\log {\bf Q}^{(\epsilon,n)}\{ 
G\}
\geq -\inf_{x \in G} \Lambda^{(n)*}(x).
 \end{equation}
\end{Th}\vspace{5mm}

\noindent {\bf Proof of Theorem (\ref{gellis})} Part 1: Upper bound.

\noindent Consider any closed and bounded set $F \subset {\bf R}^{2^n}$.

\noindent Let $\chi_A$ denote the indicator function of a set $A$. Recall that $Z_j$
are independent, identically distributed $N(0, \frac{\epsilon t}{2^n})$
random variables, so that 

\[ E[e^{pZ_j}] = \int_{-\infty}^\infty \frac{2^{n/2}}{\sqrt{2\pi \epsilon
t}}\exp\left\{-\frac{2^{n-1}x^2}{ \epsilon t} + px\right\}dx = 
\exp\left\{\frac{p^2
\epsilon t}{2^{n+1}}\right\}.\]

\noindent  Then, because the $(Z_j)_{j=1}^{2^n}$ are independent, 

\begin{eqnarray*} {\bf Q}^{(\epsilon, n)}\{ F\} & =&
E^{(\epsilon,n)}[\chi_F(Z)] \leq E^{(\epsilon,n)}\left [ \exp\left\{
\sum_{j=1}^{2^n} \lambda_j Z_j - \inf_{x \in F} \sum_{j=1}^{2^n}
\lambda_j x_j\right\}\right ] \\&  = & 
 \exp\left\{\frac{\epsilon t}{2^{n+1}}\sum_{j=1}^{2^n} \lambda_j^2 -
\inf_{x \in F}\sum_{j=1}^{2^n} \lambda_j x_j\right\}.
\end{eqnarray*}

\noindent Since $F$ is closed and bounded, there is a point $\tilde{x} \in
F$ at which the infimum is obtained. Then

\[ \epsilon \log {\bf Q}^{(\epsilon,n)}\{F\} \leq \frac{\epsilon^2
t}{2^{n+1}}\sum_{j=1}^{2^n}\lambda_j^2 - \epsilon \sum_{j=1}^{2^n}
\lambda_j\tilde{x}_j.\]

\noindent The inequality holds for all $\lambda$ and, in
particular, for $\lambda = \frac{2^n}{\epsilon t} \tilde{x}$. It follows
from equation (\ref{fenlegtran})  that for any closed bounded set $F$, 

\begin{equation}\label{cbsF} \epsilon \log {\bf Q}^{(\epsilon,n)}\{F\} \leq
-\frac{2^n}{2t}\sum_{j=1}^{2^n}\tilde{x}_j^2= - \inf_{x \in F}
\Lambda^{(n)*}(x).\end{equation}

\noindent Now the result is extended   to arbitrary closed sets. Consider a
closed set $F$ and let $\tilde{x}$ denote a point such that
$\Lambda^{(n)*}(\tilde{x}) = \inf_{x \in F}\Lambda^{(n)*}(x)$. It is easy
to see that such a point exists, because $F$ is closed and $\Lambda^{(n)*}$
is quadratic and convex. Choose a
$\delta > 0$ and choose a
$\rho$ such that
$\frac{2^n\rho^2}{2t} >
\Lambda^{(n)*}(\tilde{x}) + \delta$.  
\vspace{5mm}

\noindent  By Chebychev's inequality, for all $\alpha > 0$,

\[ {\bf Q}^{(\epsilon,n)}\{Z_j < -\rho\} =  {\bf Q}^{(\epsilon,n)}\{- Z_j >
-\rho\} \leq e^{-\alpha \rho}E \left [e^{-\alpha Z_j} \right ] = \exp\left\{-\alpha \rho +
\frac{\alpha^2}{2}\epsilon \frac{t}{2^n}\right\},\]

\noindent yielding

\[ {\bf Q}^{(\epsilon,n)}\{Z_j < -\rho\} = {\bf Q}^{(\epsilon,n)}\{Z_j >\rho\}
\leq \exp\left\{-\frac{2^n \rho^2}{2\epsilon t}\right\}\]

\noindent It follows that, for $\rho > 0$,

\[ {\bf Q}^{(n,\epsilon)}\left \{{\bf R}^{2^n} \backslash [-\rho, \rho]^{2^n} \right \}
\leq  2^{n+1} \exp\left \{-\frac{2^n \rho^2}{2\epsilon t}\right \},\]

\noindent so that 

\[ \limsup_{\epsilon \rightarrow 0} \epsilon \log {\bf Q}^{(n,\epsilon)}
\{{\bf R}^{2^n}
\backslash [-\rho,
\rho]^{2^n}\} \leq -\frac{2^n \rho^2}{2t} \leq -\inf_{x \in F}
\Lambda^{(n)*}(x) - \delta.
\]  
 
\noindent Equation (\ref{cbsF}),
  holds for all closed bounded sets $F$. This yields

\[ \limsup_{\epsilon \rightarrow 0}\epsilon \log {\bf Q}^{(n,\epsilon)}\{
F \cap [-\rho, \rho]^{2^n}\} \leq -\inf_{x \in F}\Lambda^{(n)*}(x).\] 

\noindent Choose 
$\rho$ such that
$\frac{2^n
\rho^2}{2t} >
\Lambda^{(n)*}(\tilde{x}) + \delta$. Since 
\[ \lim_{\epsilon \rightarrow
0}\epsilon \log  {\bf Q}^{(\epsilon,n)}\{ {\bf R}^{2^n}\backslash [-\rho,
\rho]^{2^n}\} \leq -\inf_{x \in F} \Lambda^{(n)*}(x),\]

\noindent it therefore follows
that  

\begin{eqnarray*} \limsup_{\epsilon \rightarrow 0}\epsilon \log
{\bf Q}^{(\epsilon,n)}\{F\} &=& 
\limsup_{\epsilon \rightarrow 0} \epsilon
\log \left( {\bf Q}^{(\epsilon,n)}\{F \cap [-\rho, \rho]^{2^n}\} +{\bf
Q}^{(\epsilon,n)}\{F\backslash [-\rho, \rho]^{2^n}\} \right)\\
&\leq &   \limsup_{\epsilon \rightarrow 0}
\epsilon
\log \left( {\bf Q}^{(\epsilon,n)}\{F \cap [-\rho, \rho]^{2^n}\} +{\bf
Q}^{(\epsilon,n)}\{{\bf R}^{2^n}\backslash [-\rho, \rho]^{2^n}\} \right)\\
 &\leq &  -\inf_{x \in F} \Lambda^{(n)*}(x).
\end{eqnarray*}

\noindent The upper bound of theorem (\ref{gellis}), given in equation
(\ref{uppbdd}), for arbitrary closed sets, has therefore been established.
\vspace{5mm}

\noindent For the lower bound advertised in equation (\ref{lowbdd}),
consider an open set
$G
\subset {\bf R}^{2^n}$. For any point $y \in G$, there exists a $\delta(y) >
0$ and a ball $B_{y,\delta(y)}$ such that $B_{y,\delta(y)} \subset G$. It is
therefore sufficient to prove that for all
$y
\in {\bf R}^{2^n}$, 

\[ \lim_{\delta \rightarrow 0}\liminf_{\epsilon \rightarrow 0}\epsilon
\log {\bf Q}^{(\epsilon, n)}\{ B_{y,\delta}\} \geq \Lambda^{(n)*}(y) =
- \frac{2^n}{2t}\sum_{j=1}^{2^n}  
 y_j^2.
\]

\noindent Let $C(n)$ denote the volume of
the unit ball in ${\bf R}^{2^n}$. Then 

\begin{eqnarray*} \lefteqn{\epsilon \log {\bf Q}^{(\epsilon,
n)}(B_{y,\delta})}\\&&=
\epsilon \log \int_{B_{y,\delta}} \frac{2^{2^{n-1} n}}{(2\pi \epsilon
t)^{2^{n-1}}
 }\exp \left\{
- \frac{2^{n-1}}{\epsilon t}\sum_{j=1}^{2^n}  z_j^2 
\right\} dz
\\ && \geq   \epsilon (2^{n-1} n\log 2 - 2^{n-1}\log (2\pi \epsilon))
  +
\epsilon
\log (C(n)
\delta^{2^n}) - \inf_{z \in B(y,\delta)}  \frac{ 2^{n-1}|z|^2}{t},
\end{eqnarray*}

\noindent yielding

\[\lim_{\delta \rightarrow 0}\left(\liminf_{\epsilon \rightarrow
0}\epsilon
\log {\bf Q}^{(\epsilon, n)}\{ B_{y,\delta}\}\right) \geq  
- \frac{2^n}{2 t}\sum_{j=1}^{2^n}   
y_j^2.
\] 

\noindent The lower bound follows directly. The proof of theorem
(\ref{gellis}) is therefore complete.
\qed \vspace{5mm}

\noindent The approach to proving theorem (\ref{varad}) is firstly to prove
a discrete version, given by proposition (\ref{laplace}) and then take a
limit, which is the subject of lemma (\ref{exchangelim}). Theorem
(\ref{varad}) then follows from proposition (\ref{laplace}) followed by
lemma (\ref{exchangelim}). 

\begin{Propn}[The Laplace Method]\label{laplace}
Let $V$ be smooth and bounded. Let $Y^{(n)}$ be defined as in equation
(\ref{wayyen}) and note that $\dot{Y}^{(n)}(s,z) = \frac{2^n}{t} z_{k(s) + 1}$,
where $\dot{Y}^{(n)}$ denotes derivative of $Y^{(n)}$ with respect to $s$.
Recall the definition of ${\cal A}$; namely, ${\cal A} : {\bf R}_+ \times
W^{1,2}({\bf R}_+) \rightarrow {\bf R}$, where 

\begin{equation}\label{action2}
 {\cal A}(t;\xi) = \frac{1}{2}\int_0^t
\dot{\xi}^2(s)ds +
\phi(\xi(0)) +
\int_0^t V(s, \xi(s))ds.\end{equation}

\noindent Note that for
$z
\in {\bf R}^{2^n}$,  
\[{\cal A} (t;Y^{(n)}(t-.,z)) =
\frac{1}{2}\sum_{j=1}^{2^n} \frac{t}{2^n}
\left( \frac{ 2^nz_j}{t}\right)^2  + \phi( Y^{(n)}(t,z)) + \int_0^tV(t-s,
Y^{(n)}(s,z))ds.
\]

\noindent Then, with $U^{(\epsilon,n)}$ defined in equation (\ref{youeeen}),

\[ \lim_{\epsilon \rightarrow 0} -\epsilon \log U^{(\epsilon,n)}(t,x) =
\inf_{z \in {\bf R}^{2^n}} {\cal A}(t; Y^{(n)}(t-.,z)).\]
\end{Propn}\vspace{5mm}

\noindent {\bf Proof}  To make the notation slightly more convenient, for $z \in {\bf R}^{2^n}$, the following is used:

\[ \tilde{\cal A}^{(n)}(z) := {\cal A}(t; Y^{(n)}(t-.,z).\]
Step 1: Upper bound. For
any open set
$G
\subset {\bf R}^{2^n}$, 

\begin{eqnarray*}  U^{(\epsilon, n)}(t,x) \geq   E^{(\epsilon,
n)}\left [e^{ -\frac{1}{\epsilon}\left(\phi(Y^{(n)}(t,Z)) +
\int_0^t V(t-s, Y^{(n)}(s,Z))ds
 \right) } \chi_G(Z_1,\ldots, Z_{2^n})\right ],
\end{eqnarray*}

\noindent from which it follows directly, by theorem (\ref{gellis})
equation (\ref{lowbdd}), that 

\begin{eqnarray*}\lefteqn{ -\liminf_{\epsilon \rightarrow 0} \epsilon \log
U^{(\epsilon, n)}(t,x)}\\&&
\leq \sup_{z \in G}\left(\phi(Y^{(n)}(t,z)) +
\int_0^t V(t-s , Y^{(n)}(s,z)) ds\right) + \inf_{z
\in G}
\frac{1}{2}\sum_{j=1}^{2^n} \frac{t}{2^n}\left(\frac{2^n
z_j}{t}\right)^2\\ 
&& \leq  \sup_{z \in G} \tilde{{\cal A}}^{(n)}(z).
\end{eqnarray*}

\noindent Set $m = \inf_{y \in {\bf R}^{2^n}} \tilde{{\cal A}}^{(n)}(y)$ and, for $\delta > 0$, 
set
$G^\delta = \{z | \tilde{{\cal A}}^{(n)}(z) - m < \delta\}$. Using the
continuity of
${\cal A}$, it is easy to see that $G^\delta$ is an open subset of ${\bf
R}^{2^n}$.  By taking $G = G^\delta$ in equation (\ref{lowbdd}) from theorem
(\ref{gellis}) and letting
$\delta \rightarrow 0$, it follows that 

\[ -\liminf\epsilon \log U^{(\epsilon,n)}(t,x) \leq \inf_{z \in {\bf R}^{2^n}}
{\cal A}(z).\]

\noindent Part 2: Lower bound. Recall that $m = \inf_{y \in {\bf
R}^{2n}}{\cal A}(y)$. Set
\[ M =
\sup_x |\phi(x)| +
\sup_{x,s} t|V(s,x)|.\]
 
\noindent Set

\[ F(z) = \phi(Y^{(n)}(t,z)) + \int_0^t  V(t-s, Y^{(n)}(s,z))ds.\] 
 Set

\[ A = \{ y | \Lambda^{(n)*}(y) \geq m+M+1\}.\]

\noindent Fix a  $\delta > 0$ and, for $0 \leq j \leq [(M+1)/\delta] + 1$,
set

\[ A^{(\delta)}_j = \{y | m + j \delta \leq \Lambda^{(n)*}(y) \leq m + (j+1)
\delta\}.\]

\noindent Set 
$N = [(M+1)/\delta] + 1$. Note that $A^{(\delta)}_0,\ldots, A^{(\delta)}_N, A$ are
closed sets and that ${\bf R}^{2^n} = \cup_{j=0}^N A^{(\delta)}_j \cup
A$. It follows that

\[ U^{(\epsilon,n)}(t,x) \leq \sum_{j=0}^N
E^{(n,\epsilon)}[e^{-\frac{1}{\epsilon} F(Z)}\chi_{A^{(\delta)}_j}(Z)] +
 E[e^{-\frac{1}{\epsilon} F(Z)}\chi_A(Z)].\]

\noindent Now, since $A_j^{(\delta)}$ is closed for each $j \in
\{0,1,\ldots, N\}$, it follows from the upper bound given by equation
(\ref{uppbdd}) in theorem (\ref{gellis})   that 

\[ -\lim_{\epsilon \rightarrow 0} \epsilon \log
E^{(n,\epsilon)}[e^{-\frac{1}{\epsilon} F(Z)}\chi_{A_{j}^{(\delta)}}(Z)] \geq
 \inf_{x \in A_j^{(\delta)}}F(x) + j \delta \geq \inf_{x \in A_j^{(\delta)}}
{\cal A}(x) - \delta
\geq m - \delta.\]

\noindent Furthermore, since $A$ is closed, it follows from the upper
bound given by equation (\ref{uppbdd}) from theorem (\ref{gellis})   that 

\[ -\liminf_{\epsilon \rightarrow +\infty}\epsilon \log
E[e^{-\frac{1}{\epsilon} F(Z)}\chi_A(Z)] \geq -M+M+m+1 = m+1.  \]

\noindent Since $A_0^{(\delta)}, \ldots, A_N^{(\delta)},A$ is a {\em
finite} collection of sets, it follows directly that 

\[ \liminf_{\epsilon \rightarrow +\infty} - \epsilon \log U^{(\epsilon,
n)}(t,x) \geq m - \delta\]

\noindent for all $\delta > 0$ and hence that 

\[ \liminf_{\epsilon \rightarrow +\infty} - \epsilon \log U^{(\epsilon,
n)}(t,x) \geq m.\]

\noindent The proof of proposition (\ref{laplace}) is complete. \qed
\vspace{5mm}

\noindent Recall that 
\begin{equation}\label{essenn} {\cal S}^{(n)} = \{ y \in C([0,t]) |  
\dot{y} = z_k,
\frac{tk}{2^n} < s < \frac{t(k+1)}{2^n}, z_k \in {\bf R} \}.
\end{equation}

\noindent Note that proposition (\ref{laplace}) has proved that 

\begin{equation}\label{reexpress} -\lim_{\epsilon \rightarrow 0} \epsilon
\log U^{(\epsilon,n)}(t,x) =
\inf_{y \in {\cal S}^{(n)}|y(t) = x}\left\{\frac{1}{2}\int_0^t
\dot{y}^2(s)ds +
\int_0^t V(s,  y(s))ds + \phi(y(0))\right\}.
\end{equation}

\noindent The aim is to let $n \rightarrow +\infty$.

\begin{Lmm}\label{exchangelim}
\[ \lim_{\epsilon \rightarrow 0} - \epsilon \log U^{(\epsilon)}(t,x) =
\lim_{n \rightarrow +\infty}\left(\lim_{\epsilon \rightarrow
0}-\epsilon \log U^{(\epsilon,n)}(t,x)\right ).\]
\end{Lmm}

\noindent {\bf Proof} Let $w^{(\epsilon)} = \sqrt{\epsilon}w$; namely, a one dimensional
Brownian motion with diffusion coefficient $\epsilon$, with
$w^{(\epsilon)}(0) \equiv 0$. For
$j=1,\ldots 2^n$, set 

\begin{equation}\label{zedhere} Z_j
= w^{(\epsilon)}\left(\frac{tj}{2^n}\right) -
w^{(\epsilon)}\left(\frac{t(j-1)}{2^n}\right).
\end{equation}

\noindent Note
that $Z_j$ are independent identically distributed random variables,
each with distribution $N(0, \frac{t\epsilon}{2^n})$. Recall the notation $Z
= (Z_j)_{j=1}^{2^n}$. Note that
$x + w^{(\epsilon)}(\frac{tk}{2^n}) = Y^{(n)}(\frac{tk}{2^n},Z)$, where
$Y^{(n)}$ is defined by equation (\ref{wayyen}). For $0 \leq s \leq t$,
set $\tilde{w}^{(\epsilon)}(s) = x + w^{(\epsilon)}(s)$. Set

\[ F(y) = \int_0^t V(t-s,y(s))ds + \phi(y(t)).\]

\noindent Then, using $E$ to denote expectation with respect to the
Brownian motion $w^{(\epsilon)}$, and using $Z$ as in equation
(\ref{zedhere}), it follows by Hölder's inequality that 

\begin{eqnarray}\nonumber \lefteqn{
\epsilon \log E\left [\exp\left\{ -\frac{1}{\epsilon}\left(\int_0^t
V(t-s, w^{(\epsilon)}(s) + x) ds + \phi ( x +
w^{(\epsilon)}(t))\right)\right\}\right ] }\\&&= \nonumber  \epsilon \log
E\left [\exp\left\{
-\frac{1}{\epsilon} F(\tilde{w}^{(\epsilon)} )\right\}\right ]\\&&
 \nonumber \leq
\frac{\epsilon n}{n+1}\log E\left [ \exp\left\{-\frac{n+1}{n\epsilon}
F(Y^{(n)}(.,Z))\right\}\right ]\\&& \label{cone} \hspace{10mm} +
\frac{\epsilon}{n+1}
\log E\left [\exp\left\{-\frac{n+1}{\epsilon}\left (
F(\tilde{w}^{(\epsilon)}) - F(Y^{(n)}(.,Z))\right)\right\}
\right ].
\end{eqnarray}

\noindent Similarly, by Hölder's inequality, it follows that  

\begin{eqnarray}\lefteqn{ \nonumber \epsilon \log E \left [ \exp
\left\{-\frac{n}{\epsilon (n+1)} F(Y^{(n)}(.,Z))\right\}\right ]  \leq
  \frac{\epsilon n}{n+1} \log E\left [
\exp\left\{-\frac{1}{\epsilon}F(\tilde{w}^{(\epsilon)})\right\}\right]}\\&&
 \label{cone2} \hspace{40mm}+
\frac{\epsilon}{n+1}\log E\left [ \exp\left\{ -
\frac{n}{\epsilon}\left (
F(Y^{(n)}(s,Z)) - F(\tilde{w}^{(\epsilon)})\right)\right\}\right ].
\end{eqnarray}

\noindent Now, set $\tilde{C} = t \sup_x \sup_s |V(s,x)| + \sup_x |\phi(x)|$ and note that

\begin{eqnarray}\nonumber \epsilon \log E \left [ \exp
\left\{-\frac{n}{\epsilon (n+1)} F(Y^{(n)}(.,Z))\right\}\right ] &\geq & \epsilon \log E \left [ \exp
\left\{-\frac{1}{\epsilon  } F(Y^{(n)}(.,Z))\right\}\right ] - \frac{\tilde{C}}{n+1}\\ & =& \epsilon \log U^{(\epsilon, n)}(t,x) - \frac{\tilde{C}}{n+1} \label{useone}
\end{eqnarray}

\noindent and 

\begin{eqnarray} \nonumber \epsilon \log E \left [ \exp
\left\{-\frac{(n+1)}{\epsilon n} F(Y^{(n)}(.,Z))\right\}\right ] &\geq & \epsilon \log E \left [ \exp
\left\{-\frac{1}{\epsilon  } F(Y^{(n)}(.,Z))\right\}\right ] + \frac{\tilde{C}}{n}\\ & =& \epsilon \log U^{(\epsilon, n)}(t,x) + \frac{\tilde{C}}{n}.\label{usetwo}
\end{eqnarray}

\noindent Set

\[ I_1 = \frac{\epsilon}{n+1}\log E \left [\exp\left\{
-\frac{n+1}{\epsilon} (F(\tilde{w}^{(\epsilon)}) - F(Y^{(n)}(.,Z)))\right\}
\right ]\]

\noindent and 

\[ I_2 = \frac{\epsilon}{n}\log E \left [\exp\left\{
-\frac{n}{\epsilon} ( F(Y^{(n)}(.,Z)) -
F(\tilde{w}^{(\epsilon)}))\right\}
\right ].\]

\noindent Note that inequalities (\ref{cone}) and (\ref{cone2}), using inequalities (\ref{useone}) and (\ref{usetwo}) now  yield

\begin{eqnarray}\lefteqn{ \nonumber \epsilon \frac{n+1}{n}\log
U^{(\epsilon,n)}(t,x) -
\frac{1}{n}\tilde{C} - I_2}\\&& \label{sandwich}\leq \epsilon \log
U^{(\epsilon)}(t,x)
\leq
\frac{\epsilon n}{n+1}\log U^{(n,\epsilon)}(t,x) + \frac{1}{n+1}\tilde{C} +
I_1.\end{eqnarray}

\noindent Now, set $C = \sup_t \sup_x |V_x(t,x)|$ and
recall that
$C < +\infty$ by the hypotheses on $V$. Set

\[ X_k = \sup_{0 \leq r \leq 1} \left | (1-r)w^{(\epsilon)}(\frac{kt}{2^n}) 
+ r w^{(\epsilon)}(\frac{(k+1)t}{2^n}) -
w^{(\epsilon)}(\frac{(k+r)t}{2^n})\right |.\]

\noindent Set 

\begin{equation}\label{waykay} \eta_k(r) :=
w^{(\epsilon)}(\frac{(k+r)t}{2^n}) -
w^{(\epsilon)}(\frac{kt}{2^n}),\end{equation}

\noindent then 

\[ X_k = \sup_{0 \leq r \leq
1}|r\eta_k(1) - \eta_k(r)|.\]

\noindent It is clear, from the basic property of Brownian
motion that increments over disjoint time intervals are independent, that   
$X_k$ are independent and identically distributed.  Since $x +
w^{(\epsilon)}(t) = Y^{(n)}(t,Z)$,  Taylor's expansion theorem, together with
the fact that
$\tilde{w}^{(\epsilon)}(\frac{kt}{2^n}) = Y^{(n)}(\frac{kt}{2^n},Z)$
for all $0 \leq k \leq 2^n$ yields  

\[
|F(Y^{(n)}(.,Z)) - F(\tilde{w}^{(\epsilon)})| \leq \frac{C}{2^n}
\sum_{k=0}^{2^n - 1}X_k,
\]

\noindent so that, for a random variable $X$ with the same distribution as
$X_k$, 

\begin{eqnarray*}\lefteqn{\frac{\epsilon}{n+1}\log E\left [ \exp\left\{ -
\frac{n+1}{\epsilon}\left (
F(Y^{(n)}(.,Z)) - F(\tilde{w}^{(\epsilon)})\right)\right\}\right ]}\\&&
\hspace{40mm} 
\leq
\frac{2^n \epsilon}{n+1}\log E\left[\exp\left\{ 
\frac{n+1}{2^n\epsilon} C X\right\}\right].
\end{eqnarray*}

\noindent Now, using $\eta_k(r)$ defined in equation (\ref{waykay}), let

\[ \tilde{\eta} =  \sup_{0 \leq r \leq 1}|\eta_k(r)|\]

\noindent and note that $2\tilde{\eta} >   X_k$.  From Revuz and
Yor~\cite{RY} page 55 proposition 1.8,

\begin{Propn} Let $\beta$ denote a standard Brownian motion with
$\beta(0) = 0$ and let $S(t) = \sup_{0 \leq s \leq t} \beta(s)$. Let ${\bf P}$ denote the probability measure under which $\beta$ is a standard Brownian motion. Then 

\[ {\bf P} (S(t) \geq a) \leq \exp \left \{ -\frac{a^2}{2t} \right \}.\]

\end{Propn}

\noindent From equation (\ref{waykay}), it follows that $\tilde{\eta}
\stackrel{(d)}{=} S\left(\frac{\epsilon t}{2^n} \right)$. Using
${\bf Q}
$ to denote the probability measure with respect to the process
$w^{(\epsilon)}$, it follows directly that 

\[ {\bf Q} (\tilde{\eta} \geq a) \leq 2\exp\left\{-\frac{2^{n-1} a^2}{\epsilon
t}\right \}.\]

\noindent Set $\gamma = \frac{2^{n-1}}{\epsilon t}$. It follows that 

\begin{eqnarray*} E\left[\exp\left\{ \alpha \tilde{\eta} \right \} \right]
&=&
\int_0^\infty {\bf Q}(e^{\alpha \tilde{\eta}} \geq x) dx \\
& = & 1 + \int_1^\infty {\bf Q}(\tilde{\eta} > \frac{\log x}{\alpha}) dx \\
&=& 1 + \int_0^\infty \alpha e^{\alpha y}{\bf Q}(\tilde{\eta} > y) dy \\
& \leq & 1 + 2\int_0^\infty \alpha e^{\alpha y} e^{-\gamma y^2} dy\\
& \leq & 1 + 2 \alpha\sqrt{\frac{\pi}{\gamma}}e^{\alpha^2 / 4\gamma}.
\end{eqnarray*} 

\noindent Using the inequality $1 + ae^b \leq e^{a+b}$ for all $a \geq 0$
and all $b \in {\bf R}$, it follows that 

\begin{eqnarray*} \frac{2^n \epsilon}{n+1}\log E[e^{\frac{n+1}{2^n \epsilon}
C X}]
  &\leq & \frac{2^n \epsilon}{n+1}\log E[e^{\frac{n+1}{2^n \epsilon}
2 C Y}]\\ & \leq & \frac{2^n \epsilon}{n+1}\log \left(1 +
\sqrt{\frac{(n+1)^2 C^2 \pi t}{2^{3n-5}\epsilon}}\exp\left
\{\frac{(n+1)^2 C^2 t}{2^{3n - 5}\epsilon}\right \}
\right)\\
& \leq & \frac{C\sqrt{\pi t \epsilon}}{2^{(n-5)/2}} + \frac{(n+1)C^2
t}{2^{2n-5}}.
\end{eqnarray*}

\noindent It follows that 

\begin{equation}\label{eyeonebound} I_1 :=   \frac{\epsilon}{n+1}
\log E\left [\exp\left\{-\frac{n+1}{\epsilon}\left (
F(\tilde{w}^{(\epsilon)}) - F(Y^{(n)}(.,Z))\right)\right\}
\right ] \leq \frac{(n+1) C^2 t}{2^{2n - 5}} + \frac{C\sqrt{\pi t
\epsilon}}{2^{(n-5)/2}}
\end{equation}

\noindent and, similarly, that 

\begin{equation}\label{eyetwobound} I_2 :=   
\frac{\epsilon}{n}
\log E\left [\exp\left\{-\frac{n}{\epsilon}\left (
F(Y^{(n)}(.,Z))- F(\tilde{w}^{(\epsilon)})\right)\right\}
\right ]   \leq \frac{n C^2 t}{2^n } + \frac{C\sqrt{\pi t
\epsilon}}{2^{(n-5)/2}}.
\end{equation}

\noindent Putting this into the inequalities (\ref{sandwich}) yields lemma
(\ref{exchangelim}) directly. \qed \vspace{5mm}

\noindent {\bf Proof of theorem (\ref{varad})}Recall that 

\begin{eqnarray*}\lefteqn{ \lim_{\epsilon \rightarrow 0}-\epsilon \log
E[e^{-\frac{1}{\epsilon}(\int_0^t V(t-s, x+w^{(\epsilon)}(s))ds + \phi(x +
w^{(\epsilon)}(t))}]}\\&& \hspace{20mm} = \lim_{n \rightarrow +\infty}
\left(\lim_{\epsilon \rightarrow 0} -\epsilon
\log E[e^{-\frac{1}{\epsilon}(\int_0^t V(t-s, Y^{(n,\epsilon)}(x,s))ds +
\phi(Y^{(n,\epsilon)}(x,t))}]\right). 
\end{eqnarray*}

\noindent Recall that 

\[{\cal S}_n = \left \{y \in C([0,t]) |  
\exists (z_1,\ldots, z_{2^n}) \in {\bf R}^{2^n}| y(s) = y(\frac{kt}{2^n}) +
(s-\frac{kt}{2^n})z_{k+1},
\frac{kt}{2^n}
\leq s \leq
\frac{(k+1)t}{2^n}\right \},\]

\noindent and that 

\begin{equation} \label{minimum}\lim_{\epsilon \rightarrow 0} - \epsilon
\log U^{(\epsilon,n)}(t,x) =
\inf_{y
\in {\cal S}_n | y(t) = x} {\cal A}(y),\end{equation}

\noindent where 

\[ {\cal A}(y) = \left\{\frac{1}{2}\int_0^t \dot{y}^2(s)ds +
\int_0^t V(t-s, y(s))ds + \phi(y(0))\right\}.\]

\noindent It follows directly from the fact that ${\cal S}_n \subset
{\cal S}_m$ for $m > n$, together with the analysis given above, that 

\[ \lim_{\epsilon \rightarrow 0} -\epsilon \log U^{(\epsilon)}(t,x) =
\inf_n \inf_{y \in {\cal S}_n} {\cal A}(y) = \inf_{y \in {\cal S}}{\cal A}(y)\]

\noindent and theorem (\ref{varad}) is proved. \qed\vspace{5mm}

\noindent It is now shown that, assuming Tychonov's
theorem, hence relative weak compactness of the unit ball in
$L^2$, the minimiser exists. It is then shown that if the
minimiser exists, then it satisfies the Euler Lagrange
equations.\vspace{5mm}

\begin{Th} [Existence of the Minimiser] \label{minimiser} Consider the
action functional

\[ {\cal A}(y) = \frac{1}{2}\int_0^t \dot{y}^2(s)ds + \int_0^t V(s,y(s))ds
+ \phi(y(0)).\]

\noindent Then, using the fact that a ball of finite radius in $L^2$ is  
compact in the relative weak topology, there exists a trajectory $\tilde{y}$
such that 
\[ {\cal A}(\tilde{y}) = \inf_{ y \in W^{1,2}([0,t])| y(t) = x}{\cal A}(y).\] 
\end{Th}

\noindent {\bf Proof of theorem (\ref{minimiser})} Consider a sequence
$(y_n)_{n=1}^\infty$ where each $y_n \in W^{1,2}([0,t])$ and $y_n(0) =
x$, such that ${\cal A} (y_1) = C < +\infty$, such that ${\cal A}(y_n)$ is
decreasing and such that 
\[ \lim_{n \rightarrow +\infty}{\cal A}(y_n) = \inf_{ y \in W^{1,2}([0,t])|
y(t) = x}{\cal A}(y).\]

\noindent Consider the sequence $(\dot{y}_n)_{n=1}^\infty$, and take a
subsequence $(\dot{y}_{n_k})_{k \geq 1}$ that is convergent to a limit $\dot{\tilde{y}}$ in the relative weak topology. That is, for any test function $g \in L^2([0,t])$, $\int_0^t g(s) \dot{y}_{n_k}(s) ds \stackrel{k \rightarrow +\infty}{\longrightarrow} \int_0^t g(s) \dot{\tilde{y}}_{n_k}(s) ds$. Such a
sequence exists, since

\begin {equation}\label{sequex}\sup_n \frac{1}{2}\int_0^t \dot{y}_n^2(s) ds \leq C + \|V\|_\infty t +
\|\phi\|_\infty < +\infty
\end{equation}

\noindent and because any ball of finite radius in $L^2$ is  compact in the relative weak topology. Let  $\tilde{y}$
denote the function such that $\tilde{y}(0) = x$, with derivative $\dot{\tilde{y}}$. It
follows by choosing test functions $\chi_{[0,s]}$, which are clearly in $L^2([0,t])$, that 

\[ \lim_{k \rightarrow +\infty} |y_{n_k}(s) - \tilde{y}(s)| =
0 \quad \forall s \in [0,t]\]

\noindent and, furthermore, using equation (\ref{sequex}) and Hölder's inequality,

\[\sup_k \sup_{0 \leq s \leq t}|y_{n_k}(s) - \tilde{y}(s)|\leq \int_0^t |\dot{y}_{n_k}(s) - \dot{\tilde{y}}(s)|ds \leq 2t^{1/2}(C + \|V\|_\infty + \|\phi\|_\infty)^{1/2}. \]

\noindent and hence, because $V$ and $\phi$ are smooth and uniformly
bounded with uniformly bounded first derivatives, it follows by the dominated convergence theorem that 

\begin{eqnarray*} \lefteqn{  \left |\left(\int_0^t V(t-s,y_{n_k}(s))ds +
\phi (y_{n_k}(t))\right) - \left(\int_0^t V(t-s,\tilde{y}(s))ds - \phi
(\tilde{y}(t))\right)
\right |}
\\&& \leq  \|V^\prime\|_\infty \int_0^t |y_{n_k} (s)- \tilde{y}(s)| ds + \|\phi_x\|_\infty |y_{n_k} (t)- \tilde{y}(t)|\\&& \stackrel{k \rightarrow +\infty}{\longrightarrow} 0.
\end{eqnarray*}

\noindent Since $\dot{y}_{n_k} \rightarrow \dot{y}$ in $L^2$ with the relative weak
topology, it follows from standard results that convergence is almost everywhere and hence, by Fatou's lemma, 

\[ \liminf_{k \rightarrow +\infty} \frac{1}{2}\int_0^t \dot{y}_{n_k}^2(s)ds
\geq \frac{1}{2}\int_0^t \dot{\tilde{y}}^2(s)ds.\]

\noindent It therefore follows that 

\[ \inf_{ y \in W^{1,2}([0,t])|
y(t) = x}{\cal A}(y) \geq {\cal A}(\tilde{y})\]

\noindent and hence that the trajectory $\tilde{y}$ is a minimiser; 

\[  {\cal A}(\tilde{y}) = \inf_{ y \in W^{1,2}([0,t])| y(t) = x}{\cal A}(y).\]

\noindent The proof of theorem (\ref{minimiser}) is complete. \qed
\vspace{5mm}

\begin{Th}\label{elsat} For $0 \leq s \leq t$, the minimiser
$\tilde{y}$ in theorem (\ref{minimiser}) satisfies the Euler Lagrange
equations

\[ \left\{\begin{array}{l} \ddot{\tilde{y}}(s) = V^\prime( s,\tilde{y}(s)) \\
\tilde{y}(t) = x,\qquad \dot{\tilde{y}}(0) =   \phi(\tilde{y}(0)).
\end{array}\right.\]
\end{Th}\vspace{5mm}

\noindent {\bf Proof of theorem (\ref{elsat})} The proof is
sketched; in this case, the fact that the minimiser satisfies the
Euler Lagrange equation is standard, following the arguments of
Tonelli in~\cite{T1,T2}. Consider any
$z
\in W^{1,2}([0,t])$ with $z(t) = 0$, so that $\tilde{y} + \epsilon z \in
W^{1,2}([0,t])$ for all
$\epsilon \geq 0$. Then

\[ \lim_{ \epsilon \rightarrow 0}\frac{{\cal A}(\tilde{y} + \epsilon z) -
{\cal A}(\tilde{y})}{\epsilon} = \int_0^t \dot{z}(s)\dot{\tilde{y}}(s) ds +
\int_0^t z(s) V_x( s,\tilde{y}(s))ds + z(0)\phi^\prime (\tilde{y}(0)).\]

\noindent Since, for any $z \in W^{1,2}([0,t])$ with $z(t) = 0$, $ \lim_{
\epsilon
\rightarrow 0}\frac{{\cal A}(\tilde{y} + \epsilon z) - {\cal
A}(\tilde{y})}{\epsilon} = 0$, it follows that for all $z \in W^{1,2}([0,t])$
with $z(t) = 0$, 

\[ 0 = z(0)\left (- \dot{\tilde{y}}(0) + \phi^\prime (\tilde{y}(0))\right) -
\int_0^t z(s)\left(\ddot{\tilde{y}}(s) - V_x( s,\tilde{y}(s))
\right) ds,
\]

\noindent from which theorem (\ref{elsat}) follows directly. \qed
\vspace{5mm}

\noindent It now remains to identify $u(t,x) = \dot{\eta}(t,\eta^{-1}(t,x))$
where $\eta$ solves equation (\ref{elnew}) and {\em minimises} the action functional.
\vspace{5mm}

\begin{Th}\label{identifyth} Let $V \in C^{\infty}({\bf R}^2)$ be smooth and
bounded and let
$\phi \in C^\infty ({\bf R})$ be smooth and bounded. Let
$u^{(\epsilon)}$ solve

\[ \left\{\begin{array}{l} u_t^{(\epsilon)} + \frac{1}{2}(u^{(\epsilon)2})_x =
\frac{\epsilon}{2}u^{(\epsilon)}_{xx} + V_x
\\ u_0^{(\epsilon)} =
\phi_x.\end{array}\right.\]

\noindent Let $u$ denote any weak in $L^2$ limit point of
$(u^{(\epsilon)})_{\epsilon > 0}$. Then there is a representation of the
solution $u$, $u(t,x) = \dot{\xi}^{(t,x)}(t)$, where $\xi$ is the 
trajectory which provides a {\em minimum} for the action

\[ {\cal A}(t,\xi) := \frac{1}{2}\int_0^t \dot{\xi}^2(s)ds + \int_0^t
V(s,\xi(s))ds +
\phi(\xi(0))\]

\noindent subject to the constraint $\xi^{(t,x)}(t) = x$, $\xi \in W^{1,2}([0,t])$.  For $0 \leq s \leq t$, this trajectory $\xi$ satisfies the Euler
Lagrange equation

\begin{equation}\label{eleqn} \left\{\begin{array}{l}\ddot{\xi} (s) =
V_x(s,\xi(s))
\\
\xi(t) = x, \qquad \dot{\xi}(0) = \phi_x(\xi(0)). \end{array}\right.
\end{equation}

\end{Th}

\noindent {\bf Proof of theorem (\ref{identifyth})} It has already been
shown that any limit point $u$ has the representation  

\[ u(t,x) = \frac{\partial}{\partial x}  
\inf_{\xi : \xi(t) = x}\left\{ \phi(\xi (0)) + \int_0^t
(\frac{1}{2}\dot{\xi}^2(s)   + V(s,\xi (s)))ds)\right \}.\]

\noindent Furthermore, it has been shown that the minimising
trajectory $\xi$ exists and satisfies the Euler Lagrange equations
(\ref{eleqn}). The identification that
$u(t,x) =
\dot{\xi}(t)$ where
$\xi$ is a {\em minimising} trajectory subject to the constraint that
$\xi(t) = x$  is completed as follows: Let $(\psi(s,t,x))_{s = 0}^t$ denote a
trajectory that minimises

\[ \phi(\xi (0)) + \int_0^t (\frac{1}{2}\dot{\xi}^2(s)   + V(s,\xi (s)))ds\]
 
\noindent subject to the conditions $\xi(t) = x$. Then, using the arguments of the proof of theorem (\ref{elsat}), the
variational calculus yields that
$\psi$ 
 satisfies 

\[ \psi_{ss}(s,t,x) =  V^\prime (s, \psi(s,t,x)),\]

\noindent  where   $V^\prime $ denotes (as usual) 
derivative of $V$ with respect to the second argument of $V$ and $\psi_s(0,t,x) =
\phi^\prime(\psi(0,t,x))$.
Furthermore,
$\psi(t,t,x)
\equiv x$, so that
$\psi_x(t,t,x)
\equiv 1$. It follows, using integration by parts, that 

\begin{eqnarray*}\lefteqn{ u(t,x) = \frac{\partial}{\partial x} \left(
\phi(\psi(0,t,x)) +
\int_0^t (\frac{1}{2}\psi_s(s,t,x)^2 + V(s,\psi(s,t,x)))ds\right )  }\\&&
= \phi^\prime (\psi(0,t,x))\psi_x(0,t,x) + [\psi_x(s,t,x)\psi_s(s,t,x)]_{s=0}^t\\&&
\hspace{30mm}  -
\int_0^t \psi_x(s,t,x)(\psi_{ss}(s,t,x) - V_x(s,t,x)) ds \\&&
= \psi_x(t,t,x)\psi_s(t,t,x) + \psi_x(0,t,x)\left( \phi^\prime (\psi(0,t,x)) -
\psi_s(0,t,x)\right) \\&&
= \psi_s(t,t,x),
\end{eqnarray*}

\noindent which is the advertised result. Theorem (\ref{identifyth}) is
proved.\qed \vspace{5mm}

\noindent One final lemma is required to finish this section, which will be
used in the sequel. In section (\ref{final}), {\em periodic} solutions to the
Burgers equation will be considered, therefore the initial condition will
depend on $\epsilon$.

\begin{Lmm}\label{samelim} Let $w^{(\epsilon)}$ denote a Brownian motion
with diffusion $\epsilon$, with $w^{(\epsilon)}(0) = 0$. Set
\[ U^{(\epsilon)}(t,x; \psi) = -\epsilon \log E\left [
\exp\left\{-\frac{1}{\epsilon} \left(
\psi(x+w^{(\epsilon)}(t)) + \int_0^t
V(t-s,x+w^{(\epsilon)}(s))ds\right)\right\}\right ].
\]

\noindent Let $(\phi^{(\epsilon)})_{\epsilon > 0}$ denote a family of
functions satisfying $\sup_x   \sup_{0 < \epsilon \leq 1}
|\phi^{(\epsilon)}(x)| < C < +\infty$ for some constant $C$, and such that 

\[ \lim_{\epsilon \rightarrow 0} \sup_x |\phi^{(\epsilon)} - \phi | =
0.\]

\noindent Then, for any $T < +\infty$, 
\[ \lim_{\epsilon \rightarrow 0} \sup_x \sup_{0 \leq t \leq T} \left
|\epsilon  \log U^{(\epsilon )}(t,x;\phi^{(\epsilon )}) - \epsilon  \log 
U^{(\epsilon )}(t,x,
\phi)\right | = 0.\]
\end{Lmm}\vspace{5mm}

\noindent {\bf Proof} This follows directly from noting that

\begin{eqnarray*}\lefteqn{\epsilon  \log U^{(\epsilon )}(t,x;
\phi^{(\epsilon )}) - \epsilon  \log U^{(\epsilon )}(t,x;\phi) }\\&& =
\epsilon  \log E_{\bf P}
\left[e^{-\frac{1}{\epsilon}(\phi^{(\epsilon )} - \phi )(x
+ w^{(\epsilon)}_t)}
\left ( \frac{e^{-\frac{1}{\epsilon}\left( \phi (x +w^{(\epsilon)}(t))
+ \int_0^t V(t-s,
x+w^{(\epsilon)}(s,t)ds)\right)}}{U^{(\epsilon)}(t,x,\phi)}\right )\right]
\end{eqnarray*}

\noindent so that

\[ |\epsilon  \log U^{(\epsilon )}(t,x;
\phi^{(\epsilon )}) - \epsilon  \log U^{(\epsilon )}(t,x;\phi)| \leq \sup_x
|\phi^{(\epsilon)}(x) - \phi(x)|.\] 

\noindent It follows directly that 

\[ \sup_{0 \leq t \leq T} \sup_x |\epsilon  \log U^{(\epsilon )}(t,x,
\phi^{(\epsilon )}) - \epsilon_n \log U^{(\epsilon )}(t,x,\phi)|  \leq
\sup_x |\phi^{(\epsilon )}(x) - \phi(x)|.\]

\noindent and lemma (\ref{samelim}) follows directly. 
\qed 

\section{ The Counter Example}\label{final}  
The following example  provides an example in which  solutions to the 
Euler Lagrange equations

\[ \left\{\begin{array}{l}
 \ddot{\eta} = V^\prime (t,\eta(t))\\
\eta(0,x) = x, \quad \dot{\eta}(0,x) = \phi^\prime(x)
\end{array}
 \right.\]

\noindent  that provide the representation $u(t,x) = \dot{\eta}(t, \eta^{-1}(t,x))$ to the 
solution of the inviscid Burgers equation 

\[ \left\{\begin{array}{l} 
 u_t + \frac{1}{2}(u^2)_x = V_x \\
u(0,.) = \phi_x 
\end{array}\right.\]

\noindent  are necessarily {\em not}  the global minimisers of the
associated action functional, for any $x \in [0,2\pi]$, for all $t > T >
0$ where, in the example given,  $T = 2\pi$.  \vspace{5mm}

\noindent The potential $V(t,x) = \cos(\sin t) - \cos (x + \sin t)$ will be used, so
 that
\[V_x(t,x) = \sin (x + \sin t)\]

\noindent and $\phi$ will be chosen as the unique initial condition so that the solution $u$ is $2\pi$ periodic in both variables. \vspace{5mm}

\noindent The viscous Burgers' equation under consideration is therefore

\begin{equation}\label{visbur}
\left\{ \begin{array}{l} u^{(\epsilon)}_t +
\frac{1}{2}(u^{(\epsilon)2})_x = \frac{\epsilon}{2}
u^{(\epsilon)}_{xx} + \sin (x + \sin (t)) \\ u^{(\epsilon)}(0,x) =
\phi^{(\epsilon)}_x (x), \end{array}\right.
\end{equation}

\noindent where $\phi^{(\epsilon)}$ will be chosen to provide space / time
periodic solutions and the Inviscid Burgers' equation under consideration  is
the viscosity limit, which satisfies   

\begin{equation}\label{eul} \left\{\begin{array}{l} u_t +
\frac{1}{2}(u^2)_x = \sin (x + \sin t) \\
u (0,x)  = \phi_x(x),  \end{array}\right. \end{equation}

\noindent where $\phi$ is the limit of $\phi^{(\epsilon)}$ and provides space / time
periodic solutions.  It will be shown later that for all $\epsilon \geq 0$, there exists a {\em unique} solution to equation (\ref{visbur}) and for all 
$\epsilon \geq 0$, there exists exactly one function
$\phi^{(\epsilon)}$ that yields periodic solutions of equation
(\ref{visbur}). The functions $\phi^{(\epsilon)}$ have a limit $\phi$, such that $\sup_{x \in [0,2\pi)}\lim_{\epsilon \rightarrow 0}|\phi^{(\epsilon)}(x) - \phi(x)| = 0$, and $\phi_x$ is the unique initial condition that yields periodic
solutions to the inviscid Burgers equation (\ref{eul}). For the periodic solutions $u^{(\epsilon)}$, there is a unique viscosity limit $u$, which is the unique periodic solution to equation (\ref{eul}). 
Attention is restricted to {\em periodic} solutions to equations
(\ref{visbur}) and (\ref{eul}). \vspace{5mm}

\noindent The associated {\em action functional} is 

\begin{equation}\label{actionex} {\cal A}(\xi; t) =  \int_0^t \left\{
\frac{1}{2}\dot{\xi}^2(s)    + \cos (\sin  s) - \cos ( \xi(s) + \sin s) \right\} ds +
\phi (\xi(0)).  \end{equation}

\noindent  Let $\phi^\prime (x)  = \phi_x(x)$. Using subscripts to
denote derivatives with respect to the subscripted variable, easy
variational calculus arguments yield that the critical points of the action
functional with constraint
$\xi(t) = x$ satisfy
$\xi(s) = \xi(s;t,x)$ for $0 \leq s \leq t$, where $\xi(.;t,x)$ satisfies

\begin{equation}\label{assoce-l}  
\left\{\begin{array}{l} \xi_{ss}(s;t,x) = \sin (\xi(s;t,x) + \sin s) \\
\xi(t;t,x) = x,\; \xi_s(0;t,x) = \phi^\prime  (\xi(0;t,x)).\end{array}
\right.
\end{equation}
 
\noindent Now suppose that equation (\ref{eul}) with initial condition
$\phi  =
\phi^{(1)}$ where that $\phi^{(1)}$ is differentiable at $0$ and
satisfies
$\phi^{(1)}_x(0) = 0$.   Then   for $t = 2n\pi$ and $x=0$,

\begin{equation}\label{fix} \xi(s;2n\pi,0) = -2n\pi + s - \sin s
\end{equation}

\noindent yields a solution to equation
(\ref{assoce-l}).  
 
\noindent Now consider equation
(\ref{eul}) with  boundary data
$\phi  =\phi^{(2)}$, where $\phi^{(2)}$ is differentiable at $\pi$ with
$\phi_x^{(2)}(\pi) = -2$.   Then, for $x=\pi$ and $t = 2n+1$, the function

\begin{equation}\label{fix2} \xi(s;(2n+1)\pi,\pi) = (2n+1)\pi   - s - \sin s 
\end{equation}

\noindent yields a solution to equation (\ref{assoce-l}). That equations
(\ref{fix}) and (\ref{fix2}) provide solutions to equation
(\ref{assoce-l}) for the prescribed boundary conditions   can be seen by
plugging into both sides. 

\begin{Lmm}\label{EZLMM} For all $n \geq 1$, $\xi(s) =
\xi(s;2n\pi;0)$, where $\xi$ is given by equation (\ref{fix}) does not
minimise the associated action functional (\ref{actionex}) with the
conditions $t = 2n\pi$, $\xi(2n\pi) = 0$ and $\phi_x(0) = 0$. For all $n \geq
1$,
$\xi(s) =
\xi(s;(2n +1)\pi, \pi)$ where $\xi$ is given by equation (\ref{fix2}) does not
provide a minimiser   for the associated action functional (\ref{actionex}),
with the conditions $t = (2n + 1)\pi$, $x = \pi$ and $\phi_x(\pi) = -2$.  
\end{Lmm}

\noindent {\bf Proof of lemma (\ref{EZLMM})} Solutions (\ref{fix})
and (\ref{fix2}) are considered separately. For solution (\ref{fix}), times $t = 2n\pi$ are chosen for
integer
$n$ and final condition $\xi (2 n \pi)  =   0$. This gives

\begin{eqnarray}\nonumber {\cal A}(\xi;0,2 \pi n) &=&
\frac{1}{2}\int_0^{2\pi n} (1 -
\cos   s)^2 ds + \int_0^{2\pi n} \cos( \sin s) ds - \int_0^{2\pi n} \cos  ( s )
ds+ \phi (0) \\ & =& \label{nonmin} \frac{3 \pi n}{2} + \int_0^{2 \pi n
}\cos(\sin   s)ds+ \phi (0)  > \frac{3 \pi n}{2} + \phi (0). 
\end{eqnarray}

\noindent  It is easy to see that the trajectory $\psi$ such that $\psi (t) \equiv 0$ $\forall t \geq 0$   is not a solution of the Euler Lagrange equation. The action is 

\begin{equation}\label{lower} {\cal A}(\psi;0,2 \pi n) = \int_0^{2 \pi n}
\cos(\sin  s)ds -
\int_0^{2
\pi n} \cos (  \sin s) ds + \phi (0)  = \phi (0),
\end{equation}

\noindent    so that ${\cal A}(\psi;0,2 \pi n) < {\cal A}(\xi;0,2 \pi n)$.The
statement in lemma (\ref{EZLMM}) connected with equation (\ref{fix}) is
now proved. Note that these two trajectories have the same boundary data
$\psi (2 \pi n) = \xi (2 \pi n) = 0$ and $\dot{\psi} (0) = \dot{\xi} (0) =
0$.\vspace{5mm}

\noindent For
equation (\ref{fix2}),  times $t = (2n+1)\pi$ are considered for integer $n$ and final
condition
$\xi (2(n+1)\pi) = 0$ is considered. The action is

\begin{eqnarray}\nonumber {\cal A}(\xi;0, (2 \pi + 1) n) &=&
\frac{1}{2}\int_0^{(2n + 1)\pi } (1 -
\cos   s)^2 ds + \int_0^{(2n + 1)\pi} \cos( \sin s) ds\\&&\nonumber  -
\int_0^{(2n + 1)\pi } \cos  ( s ) ds +\phi (\pi)\\ & =& \nonumber 
\frac{(6n + 3) \pi }{4} +
\int_0^{(2n+1)
\pi}
\cos(\sin   s)ds + \phi (\pi) \\ & > & \label{nonmin2}\frac{(6n+3) \pi
}{4} +
\phi (\pi). 
\end{eqnarray}

\noindent Note that $\dot{\xi}(0) = -2$. Compare with the trajectory 

\[\psi(t) = \left\{ \begin{array}{ll}\pi - 2t, & 0 \leq t \leq \frac{\pi}{2} \\
0, & t > \frac{\pi}{2}
\end{array}\right. \] 

\noindent Then $\psi((2n+1)\pi) = \xi ((2n+1)\pi) = 0$ and $\dot{\psi}(0)
= \dot{\xi}(0) = -2$ and, for all $n \geq 0$, 

\begin{eqnarray*}
 {\cal A}(\psi; 0,(2n+1)\pi) &=& \pi + \int_0^{\pi/2} \cos
(\sin s) ds  -
\int_0^{\pi / 2}
\cos (\pi-2s + \sin s)ds + \phi (\pi)
\\  &\leq& 2\pi + \phi (\pi).
\end{eqnarray*}

\noindent The right hand side is bounded independent of
$n$ and the action in (\ref{nonmin2}) is increasing linearly in
$n$. The statement in lemma (\ref{EZLMM})  for equation (\ref{fix2}) also
holds. The proof of lemma (\ref{EZLMM}) is complete. \qed 
\vspace{5mm}

\noindent    The crucial point is to show  that
one of the trajectories $\xi$ given by equation (\ref{fix}) or (\ref{fix2}) is
{\em necessarily} a trajectory connected with solutions of the inviscid
Burger's equation.  
\vspace{5mm}

\noindent Let $\eta(t;q,p)$ denote the solution to the equation

\begin{equation}\label{eulerequation}
\left\{\begin{array}{l}\ddot{\eta} = \sin( \eta + \sin (t)) \\ \eta(0) =
q,\qquad
\dot{\eta}(0) = p.\end{array}\right. 
\end{equation}

\noindent Several lemmas are required. The line of approach is as follows:
firstly, initial conditions $(q,p)$ which yield {\em periodic}
(modulo $2\pi$) solutions to the Euler Lagrange equations
(\ref{eulerequation}) are considered. It is shown that $(q,p) = (0,0)$ and
$(q,p) = (\pi,-2)$ are the {\em only} two, yielding solutions (modulo
$2\pi$)
$\eta(t) = t -
\sin(t)$ and
$\eta(t) = \pi - t - \sin(t)$ respectively. Next, it is shown that
there exists exactly one space / time periodic solution to the
inviscid Burgers equation and that this is obtained as the viscosity
limit of periodic solutions to the viscous Burgers equations. Finally,
it is shown that a periodic solution to the Euler Lagrange equation
is necessarily associated with the periodic solution to the inviscid
Burgers equation in the construction given by equation (\ref{rep}). It
has been shown in lemma (\ref{EZLMM}) that neither of the periodic
solutions minimise the action functional. \vspace{5mm}

\noindent Firstly, let $\eta$ solve equation (\ref{eulerequation}) and set
$Y(t) =
\eta(t) -
\sin(t) + t$ and set $X(t) = Y(-t)$. Note that $X$ satisfies

\[ \ddot{X} = \sin(X-t) + \sin(t).\]
\vspace{5mm}

\begin{Lmm}\label{ONEsol} Let $X$ denote solution to the equation

\begin{equation}\label{xback}
\left\{\begin{array}{l} \ddot{X} = \sin (X-t) + \sin(t)\\ X(0) = x, \; \dot{X}(0)
= y. \end{array}\right. 
\end{equation}

\noindent where $X : {\bf R}\rightarrow {\bf S}^1 = [0,2\pi)$; that is, $x =
x+2\pi$. Only the initial conditions $(x,y) = (0,0)$ and $(x,y) = (\pi,2)$
yield periodic solutions of period $2\pi k$ for some $k \in {\bf Z}$. The
corresponding periodic solutions on ${\bf S}^1$ are of period $2\pi$  and are

\[ X(t) \equiv 0\]
\noindent and 
\[X(t) = \pi + 2t.\]

\end{Lmm}\vspace{5mm}

\noindent {\bf Proof} Consider solutions to equation (\ref{xback}) in ${\bf
R}$. A solution of period $2\pi l$ in ${\bf S}^1$ will satisfy

\[ X(t) = \frac{C}{l}t + \sum_{k=-\infty}^\infty \alpha_k e^{i\frac{k}{l}t} \]

\noindent for some integer $C$ and some collection
$(\alpha_k)_{k=-\infty}^\infty$ such that $\alpha_k = \alpha_{-k}^*$,
where $\alpha_k^*$ is used to denote the complex conjugate of $\alpha_k$.
Set

\begin{equation}\label{aaafun} A(t,\alpha) := \sum_{k=-\infty}^\infty
\alpha_k e^{ikt}.\end{equation} 

\noindent Equation (\ref{xback}) yields

\begin{equation}\label{prefour}  -\sum_{k=-\infty}^\infty
\left(\frac{k}{l}\right)^2
\alpha_k e^{i\frac{k}{l}t} = \sin \left( \left (\frac{C}{l}-1\right )t +
\sum_{k=-\infty}^\infty \alpha_k e^{i\frac{k}{l}t}\right) + \sin
(t).\end{equation}

\noindent Using $K$ to denote the Kroneker delta function

\[ K_l(k) = \left\{ \begin{array}{ll}1 & l=k \\ 0 & l \neq k,
\end{array}\right.
\] 

\noindent it follows directly from equation (\ref{prefour}) that

\begin{equation}\label{kayy} -\left(\frac{k}{l}\right)^2 \alpha_k =
\frac{1}{2\pi}
\int_0^{2\pi} e^{-ikt}
\sin ((C-l)t + A(t,\alpha )) dt -
\frac{i}{2}(K_l(k) - K_{-l}(k)).
\end{equation}

\noindent Set 

\begin{equation}\label{phidef} \phi(\alpha) =
\frac{1}{2\pi}\int_0^{2\pi}\cos ((C-l)t + A(t,\alpha))dt.
\end{equation}

\noindent Let $(\tilde{\alpha})_{k=-\infty}^\infty$ denote a solution to
equation (\ref{kayy}). Then an easy differentiation of equation
(\ref{phidef}) yields

\begin{equation}\label{enodtog} \left(\frac{k}{l}\right)^2 \tilde{\alpha}_k =
  \frac{\partial}{\partial
\alpha_{-k}}\phi(\tilde{\alpha}) +
\frac{i}{2}(K_l(k) - K_{-l}(k)).
\end{equation}

\noindent Taking a power series expansion of $\phi$ yields 

\begin{eqnarray}\nonumber\lefteqn{ \phi(\alpha) =
\frac{1}{2}\left\{e^{i\alpha_0}\sum_{n=0}^\infty
\frac{i^n}{n!}
\sum_{k_1+ \ldots + k_n = -(C-l);k_j \neq 0} \alpha_{k_1}\ldots
\alpha_{k_n}\right. }\\&& \left. \label{phialpha} \hspace{20mm} +
e^{-i\alpha_0}\sum_{n=0}^\infty \frac{(-i)^n}{n!}\sum_{k_1+ \ldots + k_n =
(C-l);k_j \neq 0}\alpha_{k_1}\ldots \alpha_{k_n}
\right\}.
\end{eqnarray}

\noindent For $\tilde{\alpha}$ such that $X(t) = \frac{C}{l}t +
\sum_{k=-\infty}^\infty \tilde{\alpha}_k e^{i \frac{k}{l} t}$, where $X$
is a real solution to equation (\ref{xback}), it is clear (using $\ddot{X} =
\sin(X-t) + \sin(t)$) that
$|\ddot{X}|
\leq 2$, from which it follows directly that 
 
\begin{equation}\label{alphbd} |\tilde{\alpha}_k| \leq
\frac{2l^2}{k^2} \end{equation} 

\noindent for each $k \neq 0$. It  also follows easily from equation
(\ref{phidef}) that  

\begin{equation}\label{upderin} \left |\frac{\partial^n}{\partial
\alpha_{k_1}\ldots
\partial
\alpha_{k_n}}\phi(\tilde{\alpha}) \right | \leq 1.\end{equation}

\noindent Taylor's expansion theorem yields

\begin{equation}\label{taylexyield} 
\phi(\alpha) = \phi(\tilde{\alpha}) +  
\sum_{n=1}^\infty \frac{1}{n!}\sum_{k_1,\ldots, k_n}(\alpha -
\tilde{\alpha})_{k_1}\ldots (\alpha -
\tilde{\alpha})_{k_n}\frac{\partial^n}{\partial \alpha_{k_1}\ldots \partial
\alpha_{k_n}}\phi(\tilde{\alpha}).
\end{equation}

\noindent The expansion given by equation (\ref{taylexyield}) is
easily justified by considering the a priori bound on the derivatives given by
inequality (\ref{upderin}) and also the a priori bound on $\alpha$ and
$\tilde{\alpha}$ given by inequality (\ref{alphbd}). Expanding this gives

\begin{eqnarray}\label{phalfpa2} \lefteqn{\phi(\alpha) =
\phi(\tilde{\alpha}) +
\sum_{n=1}^\infty
\sum_{k_1,\ldots, k_n} \alpha_{k_1}\ldots \alpha_{k_n}
\left (\frac{1}{n!}\frac{\partial^n}{\prod_{l=1}^n \partial
\alpha_{k_l}} \phi(\tilde{\alpha})\right. }\\&& \nonumber
\hspace{10mm} +\left.
\sum_{m=1}^\infty
\frac{(-1)^m}{(n+m)!}
\sum_{j_1,\ldots, j_m}\tilde{\alpha}_{j_1}\ldots
\tilde{\alpha}_{j_m} \frac{\partial^{m+n}}{\prod_{l=1}^n \partial
\alpha_{k_l}\prod_{l=1}^m \partial \alpha_{j_l}}\phi(\tilde{\alpha})\right).
\end{eqnarray}

\noindent By comparing equations (\ref{phialpha}) and (\ref{phalfpa2}), it
follows that for all
$(k_1,\ldots, k_n): k_j \neq 0$,

\begin{eqnarray}\nonumber\lefteqn{  \frac{\partial^n}{\prod_{l=1}^n
\partial
\alpha_{k_l}}\phi(\tilde{\alpha}) + \sum_{m=1}^\infty\frac{(-1)^m
n!}{(n+m)!}
\sum_{j_1,\ldots, j_m}\tilde{\alpha}_{j_1}\ldots \tilde{\alpha}_{j_m}
\frac{\partial^{m+n}}{\prod_{l=1}^n \partial
\alpha_{k_l}\prod_{l=1}^m \partial \alpha_{j_l}}\phi(\tilde{\alpha})}\\&&
\hspace{55mm}  =
\left\{\begin{array}{ll}\frac{i^n}{2} & k_1+\ldots + k_n
= -(C-l) \\ \frac{(-i)^n}{2} & k_1+\ldots + k_n = C-l \\
0 & \mbox{other}\; k_1,\ldots, k_n.\end{array}\right.\label{goodbound}
\end{eqnarray}

\noindent From equation (\ref{phidef}), note that 

\begin{equation}\label{enodd} \frac{\partial^{2n+1}}{\partial
\alpha_{k_1}\ldots
\partial
\alpha_{k_{2n+1}}} \phi(\alpha) = (-1)^n \frac{\partial}{\partial
\alpha_{\sum_{l=1}^{2n+1}k_l}}\phi(\alpha).
\end{equation}

\noindent Let $K = (2\pi +
4 l^2\sum_{j=1}^\infty \frac{1}{j^2})$, so that, using bounds (\ref{alphbd})
and (\ref{upderin}), together with $|\alpha_0| < 2\pi$,   

\[ \Theta(n) := \left |\sum_{m=1}^\infty\frac{(-1)^m n!}{(n+m)!}
\sum_{j_1,\ldots, j_m}\tilde{\alpha}_{j_1}\ldots \tilde{\alpha}_{j_m}
\frac{\partial^{m+n}}{\prod_{l=1}^n \partial
\alpha_{k_l}\prod_{l=1}^m \partial \alpha_{j_l}}\phi(\tilde{\alpha})\right | 
\leq \sum_{m=1}^\infty \frac{n!}{(n+m)!}K^m.
\]

\noindent Stirling's formula, found in~\cite{jj} page 327 yields that 

\[ \frac{n!}{(n+m)!} = e^m \left(\frac{n}{n+m}\right)^{n+\frac{1}{2}}
\frac{1}{(n+m)^{m+\frac{1}{2}}}e^{\theta_{m,n}},\]

\noindent where $-\frac{1}{12(n+m)} \leq \theta_{m,n} \leq \frac{1}{12n}$,
so that 

\[ \Theta(n) \leq e\sum_{m=1}^\infty \frac{(Ke)^m}{(n+m)^{m+\frac{1}{2}}}
\stackrel{n \rightarrow +\infty}{\longrightarrow} 0.\]

\noindent It follows from equation (\ref{goodbound}), together with
equation (\ref{enodd}), that  for {\em all} $k \in {\bf Z}$, 

\[\frac{\partial}{\partial \alpha_k}\phi(\tilde{\alpha}) =
\left\{\begin{array}{ll} \frac{i}{2}e^{i\alpha_0} & k = -(C-l) \\
-\frac{i}{2}e^{-i\alpha_0} & k = C-l \\ 0 & k \neq \pm (C-l).
\end{array}\right.\]

\noindent Now consider equation (\ref{enodtog}). For $C = l$ there is no
solution; $\frac{\partial}{\partial \alpha_0}\phi(\tilde{\alpha}) =
\frac{i}{2}e^{ia_0}$ and equation (\ref{enodtog}) yields (for $k = 0$)
$0 = \frac{i}{2}e^{ia_0}$. \vspace{5mm}

\noindent For
$C=0$, this yields
$\tilde{\alpha}_k = 0$ for all $k \neq 0, \pm l$ and

\[ X(t) = \alpha_0 -\frac{i}{2}(e^{i\alpha_0} - 1) e^{it} +
\frac{i}{2}(e^{-i\alpha_0} - 1)e^{-it}.\]

\noindent For this to satisfy

\[ \ddot{X} = \sin(X-t) + \sin (t)\]

\noindent requires

\[ -\sin(\alpha_0 + t) + \sin (t) = \sin(\alpha_0 - t + \sin(\alpha_0 + t) -
\sin(t)) + \sin(t),\] 

\noindent and the only solution
occurs when
$\alpha_0 = 2n\pi$, yielding $X(t) \equiv 2n\pi$ for $n \in
{\bf Z}$, from which it follows that $X(t) \equiv 0$ on ${\bf
S}^1$.\vspace{5mm}

\noindent For $C = 2l$, this yields $\tilde{\alpha}_k = 0$ for all $k \neq 0,
\pm l$ and the equations yield

\[ X(t) = \alpha_0 + 2t + \frac{i}{2}(1 + e^{i\alpha_0})e^{it} -
\frac{i}{2}(1+e^{-i\alpha_0})e^{-it}.\]

\noindent It follows that  $\alpha_0$ must satisfy

\[ \sin(\alpha_0 + t) + \sin(t) = \sin(\alpha_0 + t - \sin(t) - \sin
(\alpha_0 + t)) + \sin(t).\]

\noindent This requires $\alpha_0 = (2n+1)\pi$,   yielding

\[ X(t) = (2n+1)\pi + 2t.\]

\noindent For $C \neq 0, l, 2l$, the equations yield

\[ X(t) = \alpha_0 + \frac{C}{l}t - \sin (t) 
- \left(\frac{l}{C-l}\right)^2\sin (\alpha_0 + \frac{l}{C}t).
\]

\noindent and it is easy to see that there are no solutions to

\begin{eqnarray*}\lefteqn{ \sin(t) + \left( \frac{l}{C}\right)^2
\left(\frac{l}{C-l}\right)^2\sin (\alpha_0 + \frac{l}{C}t) }\\&& = \sin \left(
\alpha_0 + \frac{C-l}{l}t - \sin (t) 
- \left(\frac{l}{C-l}\right)^2\sin (\alpha_0 + \frac{l}{C}t)\right) +\sin(t)
\end{eqnarray*}
 for $C \neq
0,l,2l$.  The result is established. 
 \qed \vspace{5mm}

\noindent Several lemmas connected with the Burgers' equation are
necessary. They are all standard and inserted for completeness. Their aim
is to show that, with the potential under consideration, there exists a {\em
periodic} solution to the inviscid Burgers equation which arises
as a viscosity limit and that there necessarily exists one periodic solution
to the Euler Lagrange equations used in the construction of any periodic
viscosity solution to the inviscid Burgers equation. 
\vspace{5mm}

\noindent Firstly, it is established that $\int_0^{2\pi}
u^{(\epsilon)2}(t,x) dx$ is bounded, where the bound does not depend on $t$
or $\epsilon > 0$. This is standard; the proof is included for
completeness.\vspace{5mm}

\begin{Lmm} \label{ltwobound} Let $V \in C^\infty ({\bf R}_+ \times {\bf
R})$ be uniformly bounded and $2\pi$ periodic in both variables, satisfying
$\int_0^{2\pi}V(t,x) dx = 0$ for all $t \geq 0$.  Let
$u^{(\epsilon)}$ denote the solution to the equation
\begin{equation}\label{youepspre}
\left\{\begin{array}{l}u^{(\epsilon)}_t + \frac{1}{2}(u^{(\epsilon)2})_x =
\frac{\epsilon}{2}u^{(\epsilon)}_{xx} + V_x(t,x)\\
u^{(\epsilon)}(0,.) = \mbox{initial condition.}\end{array}\right. 
\end{equation}

\noindent Let $K_1 = \sup_{(s,x)}|V(s,x)|$, $K_2 = \sup_{(s,x)}|V_x(s,x)|$ and
let
$\|f\| :=
\left (\frac{1}{2\pi}
\int_0^{2\pi} f^2(x)dx
\right )^{1/2}$. If
$\|u^{(\epsilon)}(0)\| \leq 9\pi $, then $\|u^{(\epsilon)}(t)\| \leq
 6\pi K_2 + \sqrt{6 K_1 + 6 \pi K_2}$ for all
$t
\geq 0$.
\end{Lmm}\vspace{5mm}

\noindent {\bf Proof} Set $C(t) = \|u^{(\epsilon)}(t)\|^2$. Let
$v^{(\epsilon)}$ solve

\[ v^{(\epsilon)}_t = \frac{\epsilon}{2}v^{(\epsilon)}_{xx} -
\frac{1}{2}(v_x^{(\epsilon)})^2 + \frac{1}{2}C(t) - V(t,x),\]

\noindent with initial condition $v^{(\epsilon)}(0,.)$ satisfying $\int_0^{2\pi}
v^{(\epsilon)}(0,x)dx = 0$ and $v^{(\epsilon)}_x (0,x) = u^{(\epsilon)}(0,x)$.
Note that $v^{(\epsilon)}_x = u^{(\epsilon)}$ and that $\int_0^{2\pi}
v^{(\epsilon)}(t,x) dx = 0$ for all $t \geq 0$. Furthermore, note that 
$v^{(\epsilon)} = -\epsilon \log U^{(\epsilon)}$ where $U^{(\epsilon)}$
satisfies

\begin{equation}\label{bigU} \left\{ \begin{array}{l} U_t^{(\epsilon)} =
\frac{\epsilon}{2}U^{(\epsilon)}_{xx} - \frac{1}{\epsilon}
U^{(\epsilon)}\left\{ \frac{1}{2}C(t)- V(t,x)) 
\right\}\\ U^{(\epsilon)}(0,x) =
\exp\{-\frac{1}{\epsilon}v^{(\epsilon)}(0,x)\}.\end{array}\right.
\end{equation}

\noindent Let $w$ denote a standard Brownian motion, with $w_0 =
0$. Let $w_{s,t} = w_t - w_s$. Let ${\bf P}$ denote the probability
measure associated with $w$ and let $E_{\bf P}[.]$ denote expectation
with respect to ${\bf P}$. The solution to equation (\ref{bigU}) has
Kacs representation

\begin{eqnarray}\label{kacsrep} \lefteqn{
U^{(\epsilon)}(t,x) =}\\&& \nonumber E_{\bf P} \left[
\exp\left\{-\frac{1}{\epsilon}v^{(\epsilon)}(0,x+\sqrt{\epsilon}w_{0,t})
+
\frac{1}{\epsilon}\int_0^tV(s,
x+\sqrt{\epsilon}w_{s,t} )ds\right\}\right]e^{-\frac{1}{2\epsilon}
\int_0^t
C(s)ds}
\end{eqnarray}

\noindent Using representation (\ref{kacsrep}) of the solution to equation
(\ref{bigU}), it follows that 

\begin{eqnarray}\nonumber \lefteqn{ v^{(\epsilon)}(t+s,x) = 
\frac{1}{2}\int_s^{t+s} C(r) dr}\\&& \nonumber -\epsilon
\log E\left [\exp\left\{ -\frac{1}{\epsilon}\left(
 v^{(\epsilon)}(s, x+ \sqrt{\epsilon} w_{s,t+s}) - \int_s^{t+s}
V(r, x + \sqrt{\epsilon} w_{r,t+s}  )dr\right)\right\}\right ]\\&&
\label{vplugin} \geq    \frac{1}{2}\int_s^{t+s} C(r)
dr - K_1 t  -\epsilon
\log E\left [\exp\left\{ -\frac{1}{\epsilon} 
 v^{(\epsilon)}(s, x+ \sqrt{\epsilon} w_{s,t+s})  \right\}\right ]. 
\end{eqnarray}

\noindent Now, since $\int_0^{2\pi}v^{(\epsilon)}(t,x)dx = 0$ and
$v^{(\epsilon)}(t,.)$ is continuous, it follows that there exists a point $x(t)$
such that
$v^{(\epsilon)}(t,x(t)) = 0$. It follows that 

\begin{eqnarray*} \sup_x|v^{(\epsilon)}(t,x)| &=& \sup_x |\int_{x(t)}^x
u^{(\epsilon)}(t,y) dy| \\
& \leq & \int_0^{2\pi} |u^{(\epsilon)}(t,y)|dy \\
& \leq & 2\pi \|u^{(\epsilon)}(t,.)\|.
\end{eqnarray*}

\noindent Therefore, $\sup_x |v(s,x)|
\leq 2\pi 
\|u(s)\|$. Using this, together with 
$\int_0^{2\pi} v^{(\epsilon)}(t,x) dx = 0$, it follows from inequality
(\ref{vplugin}) that 

\[ 0 \geq \frac{1}{2}\int_s^{t+s} C(r)dr - K_1 t - 2\pi
\|u^{(\epsilon)}(t,.)\|.\]

\noindent This may be rewritten as 

\[ 0 \geq  \frac{1}{2}\int_s^{t+s} C(r)dr - K_1 t - 2\pi C^{1/2}(t),\]

\noindent so that 

\begin{equation}\label{grfirst} \int_s^{t+s}C(r) dr \leq 2K_1t + 4\pi 
C(s)^{1/2}.
 \end{equation}

\noindent Equation (\ref{youepspre}) yields 

\begin{eqnarray*} \frac{d}{dt} \| u(t)\|^2 &=&
-\frac{\epsilon}{2\pi}\int_0^{2\pi} |u^2_x(t,x)| dx +
\frac{1}{\pi}\int_0^{2\pi} u(t,x)V_x (t, x ) dx \\ 
& \leq & 2K_2 \| u(t)\|,
\end{eqnarray*}

\noindent so that 

\begin{equation}\label{grsecond}
 \frac{d}{dt} C(t)^{1/2}  \leq K_2.
\end{equation}

\noindent Directly from equation (\ref{grsecond}), it follows that 
$C^{1/2}(t+s)
\leq C^{1/2}(s) + K_2t$. This may be written, for
$r < t+s$, as 

\[ C^{1/2}(r) \geq C^{1/2}(t+s) - K_2 (t+s-r).\]

\noindent For $r > (t+s) - \frac{C^{1/2}(t+s)}{K_2}$, squaring both sides
yields  

\[ C(r) \geq C(t+s) - 2K_2 (t+s-r)C^{1/2}(t+s) + K_2^2 (t+s-r)^2\]

\noindent so, for $s$ such that $C^{1/2}(t+s) > K_2 t$, integration yields 

\[ \int_s^{t+s}C(r) dr \geq tC(t+s) - K_2 t^2 C^{1/2}(t+s) +
K_2^2 \frac{t^3}{3},\]

\noindent for all $t \in [0, \frac{C^{1/2}(t+s)}{K_2}]$.  The inequality
(\ref{grfirst}) may then be applied, giving 

\[  tC(t+s) - K_2 t^2 C(t+s)^{1/2}  + K_2^2\frac{t^3}{3} \leq 2K_1 t + 4\pi 
C(s)\]

\noindent for all $t \in [0,  \frac{C^{1/2}(t+s)}{K_2}]$.  For any fixed $T <
+\infty$, set
$M(T) =
\sup_{0
\leq s
\leq T} C(s)$. Choose $\tau \in [0,T]$ such that $C(\tau) = M(T)$.  then for
all
$t
\in [0, \frac{M^{1/2}(\tau)}{K_2}]$, it follows (taking $\tau$ in the place of $t+s$) that

\[ t M(\tau) - K_2 M^{1/2}(\tau) t^2 + K_2^2 \frac{t^3}{3} \leq 2K_1 t + 4\pi  M(\tau).\]

\noindent  Now choose $t = \frac{M^{1/2}(\tau)}{K_2}$, so that 

\[  M(\tau) \leq 6K_1+ 12\pi K_2 M^{1/2}(\tau)\]

\noindent yielding

\[ M^{1/2}(\tau) \leq 6\pi K_2 + \sqrt{6 K_1 + 6 \pi K_2}.\]

\noindent Since this holds for all $T$, the result  follows
directly. \qed \vspace{5mm}

\noindent The a priori bound in lemma (\ref{ltwobound}) is a useful 
step for establishing existence of {\em periodic} solutions to the Burgers'
equation under consideration.\vspace{5mm}

\begin{Lmm} Let $V$ be  smooth and bounded and periodic in the space variable and let $u_0$ be a bounded initial condition such that $\int_0^{2\pi} u_0(x)dx = 0$. Then, for each $\epsilon \geq 0$, there exists a unique solution to the equation
 
\[ \left\{\begin{array}{l}u_t^{(\epsilon)} + \frac{1}{2}(u^{(\epsilon)2})_x = \frac{\epsilon}{2}u^{(\epsilon)}_{xx} + V_x \\ u_0 = \mbox{initial condition} \end{array}
 \right.\]
\end{Lmm}

\noindent {\bf Proof} Let $u^{(1)}$ and $u^{(2)}$ denote two solutions, set $S = u^{(1)} + u^{(2)}$ and $D = u^{(1)} - u^{(2)}$. Let $D(t,x) = \sum_n \lambda_n(t) e^{inx}$. Since $\int_0^{2\pi}u_0(x)dx = 0$, it follows that $\lambda_0 \equiv 0$. Let $\tilde{D}(t,x) = -i \sum_n \frac{\lambda_n(t)}{n}e^{inx}$. Then $D = \tilde{D}_x$ and 

\[ \left\{\begin{array}{l}\tilde{D}_t^{(\epsilon)} + \frac{1}{2}S\tilde{D}_x  = \frac{\epsilon}{2}\tilde{D}_{xx} \\ \tilde{D}(0,x) \equiv 0.\end{array}
 \right.\]

\noindent Let $w^{(\epsilon)}$ denote Brownian motion with diffusion $\epsilon$ and set

\[ X_{s,t}(x) = x + (w_t^{(\epsilon)} - w_s^{(\epsilon)}) - \frac{1}{2}\int_s^t S(r,X_{r,t}(x))dr.\]

\noindent Then, let $E^{(P)}$ denote expectation with respect to $w^{(\epsilon)}$,

\[ \tilde{D}(t,x) = E^{(P)}[\tilde{D}(0, X_{0,t}(x))] \equiv 0\quad \forall t \geq 0, \quad x \in [0,2\pi)\]

\noindent so that $D \equiv 0$. This holds for all $\epsilon \geq 0$. \qed

\begin{Lmm}\label{krep} Let $S$ be a bounded function, $2\pi$ periodic in both variables, satisfying $\int_0^{2\pi} S(t,x) dx = 0$ for all $t \in [0,2\pi]$, such that $0 < \int_0^{2\pi}\int_0^{2\pi} |S(t,x)|^2 dt dx < +\infty$. Let $Z$ satisfy

\[ \left \{ \begin{array}{ll}
 \frac{\partial }{\partial s}Z_{s,t}(x) = S(s,Z_{s,t}(x)), &  (s,t) \in {\bf R}^2 \\Z_{t,t}(x) = x& t \in {\bf R}.
\end{array}
\right. \]

\noindent Then there exists a number $k \in {\bf N}$ such that for all $s \in {\bf R}$ and all $(x,y) \in [0,2\pi)^2$

\[ \lim_{t \rightarrow +\infty}|k(Z_{s,t}(x)- Z_{s,t}(y))\mod (2\pi)|= 0\]

\noindent and

\[ \lim_{t \rightarrow -\infty}|k(Z_{s,t}(x)- Z_{s,t}(y))\mod (2\pi)|= 0.\]

 \end{Lmm}\vspace{5mm}

\noindent {\bf Proof} Set

\begin{equation}\label{queue} Q_{nm}(s,s+t) = \frac{1}{2\pi}\int_0^{2\pi}\int_0^{2\pi}e^{i(mx - ny)}\delta_y (Z_{s,s+t}(x))dy dx = \frac{1}{2\pi} \int_0^{2\pi}e^{i(mx - n Z_{s,s+t}(x))}  dx.\end{equation}

\noindent Note that 

\[ Q_{0,m}(s,s+t)= \frac{1}{2\pi}\int_0^{2\pi}e^{imx}dx = 
\left \{ \begin{array}{ll} 1 & m = 0 \\ 0 & m \neq 0
\end{array}
\right.\]

\noindent It follows directly from equation (\ref{queue}) that

\begin{equation}\label{fourt}
 e^{-in Z_{s,s+t}(x)} = \sum_m Q_{nm}(s,s+t)e^{-imx},\end{equation}

\noindent from which, for all $n \in {\bf Z}$,

\[ \sum_m |Q_{nm}(s,s+t)|^2 = \frac{1}{2\pi}\int_0^{2\pi}|e^{-in Z_{s,s+t}(x)}|^2 dx = 1.\]

\noindent Furthermore, $Q_{nm}(s,s+t) = Q_{nm}(s+2\pi,s+2\pi+t)$ for all $(s,t) \in {\bf R}^2$. For any sequence $t_1,t_2,\ldots t_k$, with $t_0 = 0$ (where multiplication is taken in the sense of matrix multiplication), 

\[Q(s,s+t_1+ \ldots + t_k) = \prod_{j=0}^{k-1} Q(s + \sum_{l=0}^j t_l ,s + \sum_{l=0}^{j+1} t_l). \]

\noindent That is,

\[ Q_{nm}(s,s+\sum_{j=1}^k t_j) = \sum_{p_1,\ldots, p_k} Q_{np_1}(s,s+t_1)Q_{p_1p_2}(s+t_1,s+t_2)\ldots Q_{p_km}(s+\sum_{j=1}^{k-1}t_j, s+\sum_{j=1}^{k}t_j).\]

 \noindent For fixed $s$, let $Q$ denote $Q(s,s+2\pi)$. It is clear that one may construct a decomposition $Q = AP$, where $P$ is orthonormal, $P^*_{nm} = P_{-n,-m}$, $A_{nm} = 0$ for $|m| \geq |n| + 1$ and $A^*_{nm} = A_{-n,-m}$.  Since $\sum_m |Q_{nm}|^2 = 1$ and $P$ is orthonormal, it follows that for each $n \in {\bf Z}$, 
\[ \sum_m |A_{nm}|^2 = \sum_m |Q_{nm}|^2 = 1.\]

\noindent Note that $A_{00} = 1$ and $P_{00} = 1$.

 \noindent   Set

\[ {\cal S} = \{ n \in {\bf Z} | Q_{n0} = 0\}\]

\noindent and set ${\cal R} = {\bf Z} \backslash {\cal S}$. Set 

\[ {\cal T}^0 = \{n \in {\cal S} | Q_{nm} = 0 \; \forall m \in {\cal R}\}\]

\noindent and, for $n \geq 0$,

\[ {\cal T}^{(n+1)} = \{n \in {\cal S} | Q_{nm} = 0 \; \forall m \in {\bf Z} \backslash {\cal T}^{(n)}\}.\]

\noindent It is clear that for all $n \in {\bf N}$, ${\cal T}^{(n+1)} \subseteq {\cal T}^{(n)}$. Set $\hat{\cal S} = \cap_{n \geq 1} {\cal T}^{(n)}$.  Note that $Q^N_{n0} = 0$ for all $n \in \hat{\cal S}$ and all $N \geq 1$. Let $Q^{(S)}$ denote $\{ Q_{mn} |(m,n) \in \hat{\cal S} \times \hat{\cal S} \}$  and note that $Q^N_{nm} = Q^{(S)N}_{mn}$ for all $(m,n) \in \hat{\cal S} \times \hat{\cal S}$. Let $\tilde{\cal R} = {\bf Z} \backslash \hat{\cal S}$ and let $Q^{(R)}$ denote $\{ Q_{nm} |(m,n) \in  \tilde{\cal R} \times \tilde{\cal R}\}$ and note that $Q^N_{nm} = Q^{(R)N}_{nm}$ for $(m,n) \in \tilde{\cal R}\times \tilde{\cal R}$.  

\noindent Let $Q^{(R)N} = A_N^{(R)} P_N^{(R)}$ denote a decomposition where $P_N^{(R)}$ is orthonormal, $A_{N;nm}^{(R)} = 0$ for $|m| \geq |n| + 1$, $P_{N;nm}^{(R)*} = P_{N;-n, -m}^{(R)}$ and $A_{N;nm}^{(R)*} = A_{N;-n, -m}^{(R)}$. Note that 

\[ Q^{(R)N} = A_N^{(R)} P_N^{(R)} = (A^{(R)} P^{(R)})^N.\]

\noindent Set $\tilde{A}^{(R)} = A^{(R)-1}A_2^{(R)}$ and note that $\tilde{A}^{(R)}_{nm} = 0$ for $|m| \geq |n| + 1$. Set $\tilde{P}^{(R)} = P_2^{(R)}P^{(R)-1}$, so that $\tilde{P}^{(R)}$ is (clearly) orthonormal. Then, since  

\[ A^{(R)}P^{(R)}A^{(R)}P^{(R)} = A_2^{(R)}P_2^{(R)},\]

\noindent it follows that 

\[ P^{(R)}A^{(R)} = \tilde{A}^{(R)} \tilde{P}^{(R)}.\]

\noindent From this, it follows that 

\[ Q^{(R)2N} = (A^{(R)}P^{(R)})^{2N} = (A^{(R)}\tilde{A}^{(R)})^N (P^{(R)}\tilde{P}^{(R)})^N = A_2^{(R)N} P_2^{(R)N}.\]  

\noindent By construction,  it is clear that for each $n \in \tilde{\cal R}$, $\sum_{m= -|n|+1}^{|n| - 1} |A_{2;nm}^{(R)}|^2 > 0$, from which it follows that $\lim_{N \rightarrow +\infty} A_{2;nm}^{(R)N} = 0$ for all $m \neq 0$ and $\lim_{n \rightarrow +\infty}|A_{2;n0}^{(R)N}| = 1$ for all $n \in \tilde{\cal R}$.  By construction, $Q_{n0}^2 = A_{2; n0}$. It follows that, for all $n \in \tilde{\cal S}$, $Q^{2N}_{n0} = 0$ for all $N \geq 1$, and for all $n \in \tilde{\cal R}$, $\lim_{N \rightarrow +\infty} |Q^{2N}_{n0}| = 1$. Choose $k_1 = \inf\{ n \geq 1 | n \in \tilde{\cal R}\}$. Since

\[ e^{- ik_1 Z_{s,s+2N}(x)} = \sum_m Q_{km}e^{- imx},\]

\noindent it follows that for each $x \in [0,2\pi)$

\[ \lim_{N \rightarrow +\infty}|e^{- ik_1 Z_{s,s+4N\pi}(x)} - Q^{2N}_{n0}| = 0\]

\noindent and hence that for all $x,y$

\[ \lim_{N \rightarrow +\infty}| k_1 (Z_{s,s+4N\pi}(x) - Z_{s,s+4N\pi}(y))\mod(2\pi)| = 0,\]

\noindent from which it is easy to show that 

\[ \lim_{t \rightarrow +\infty}| k_1 (Z_{s,s+t}(x) - Z_{s,s+t}(y))\mod(2\pi)| = 0.\]

\noindent Since the same arguments work running time `backwards', there exists a $k_2 \geq 1$ such that 

\[ \lim_{t \rightarrow -\infty}| k_2(Z_{s,s+t}(x) - Z_{s,s+t}(y))\mod(2\pi)| = 0.\]

\noindent Take $K = k_1\times k_2$ and the result follows. Now suppose that $\tilde{\cal R} = \{ 0\}$.  Then, $Q_{n0} (s,s+t) = 0$ for all $t \geq 0$ $n \neq 0$. Note that 

\[ Q_{n0}(s,s+t) = \frac{1}{2\pi}\int_0^{2\pi} e^{inZ_{s,s+t}(x)}dx = \frac{1}{2\pi}\int_0^{2\pi} e^{in y} \frac{\partial Z^{-1}_{s,s+t}(y)}{\partial y}dy,\]

\noindent so that 

\[ \frac{\partial Z^{-1}_{s,s+t}(y)}{\partial y} \equiv 1.\]

\noindent Since this holds for all $t \in {\bf R}$, it follows that 

\[ Z_{s,t}(y) = c(s,t) + y.\]

\noindent It follows that 

\[ \frac{\partial}{\partial s} c(s,t) = S(s, c(s,t) + y).\]

\noindent Since $\int_0^{2\pi} S(s,y) dy = 0$ and $S$ is $2\pi$ periodic, it follows that $\frac{\partial}{\partial s}c(s,t) = 0$ and hence that $S \equiv 0$. It follows that, if the hypotheses are satisfied, then there exists a $K \in {\bf N}$ such that  

\[\lim_{t \rightarrow +\infty}|K(Z_{s,t}(x)- Z_{s,t}(y))\mod (2\pi)|= 0\] 
\noindent and
\[\lim_{t \rightarrow -\infty}|K(Z_{s,t}(x)- Z_{s,t}(y))\mod (2\pi)|= 0.\]
\noindent The result follows.
\qed 

\begin{Lmm}\label{periodexiuni}
For each fixed $\epsilon \geq 0$, there exists an initial condition $u^{(\epsilon)}(0,.)$ which provides a  solution, $2\pi$ {\em periodic} in
both
$x$ and $t$, to the equation 
\begin{equation}\label{youeps}
\left\{\begin{array}{l}u^{(\epsilon)}_t + \frac{1}{2}(u^{(\epsilon)2})_x =
\frac{\epsilon}{2}u^{(\epsilon)}_{xx} + \sin(x + \sin(t))\\
u^{(\epsilon)}(0,.) = \mbox{initial condition.}\end{array}\right. 
\end{equation}

\noindent satisfying
$\|u^{(\epsilon)}(t)\|
\leq 9\pi $ for all $t \geq 0$ and $\int_0^{2\pi} u^{(\epsilon)}(t,x) dx =
0$ for all
$t
\geq 0$. For all $\epsilon \geq 0$, there is {\em exactly one} space
time periodic solution to equation (\ref{youeps}) satisfying
 $\int_0^{2\pi} u^{(\epsilon)}(t,x) \equiv 0$. 
\end{Lmm}
\vspace{5mm}

\noindent {\bf Proof} Note that equation (\ref{youeps}) is
simply equation (\ref{youepspre}) with $V(t,x) =   -
\cos(x +
\sin(t))$. Here $K_1 = K_2 = 1$, so that the $L^2$ norm is uniformly bounded. The line of proof is as follows:
firstly, existence of a periodic solution for
$\epsilon > 0$ is established, then uniqueness for $\epsilon > 0$. With the
uniform bounds on the
$L^2$ norm, existence of periodic solution for $\epsilon = 0$ follows from
existence of periodic solution of $u^{(\epsilon)}$, together with the
uniform bounds in $\epsilon$ on the $L^2$ norm. Finally, uniqueness for
$\epsilon = 0$ is then proved. From the uniqueness results, it therefore 
follows that
$u^{(\epsilon)}$ converges in
$L^2$, in the relative weak toplogy, to $u$.
\vspace{5mm}

\noindent {\bf Part 1: Existence for $\epsilon > 0$} Let
$U^{(\epsilon)}$ satisfy the equation

\begin{equation}\label{bigU2}
 \left\{\begin{array}{l}U^{(\epsilon)}_t =
\frac{\epsilon}{2}U^{(\epsilon)}_{xx} +
\frac{1}{\epsilon}U^{(\epsilon)}\cos (x + \sin t)\\
U^{(\epsilon)}(0,.,U_0) = U_0\end{array}\right.
\end{equation}

\noindent where $U_0$ is a non negative, bounded initial condition.

\noindent Consider the operator $T : L^2({\bf S}^1) \rightarrow L^2({\bf S}^1)$ defined by

\[ T(\phi)(x) := E\left [ \phi (x + w^{(\epsilon)}_{2\pi}) \exp\left \{\frac{1}{\epsilon}\int_0^{2\pi} \cos (x + w^{(\epsilon)}_{2\pi - s} + \sin(s))ds\right \}\right ],\] 

\noindent where ${\bf S}^1$  here (as stated before) is the circle $[0,2\pi)$; the real line with the identification $x = x + 2\pi$. Note that 

\[ T(\phi)  = U^{(\epsilon)}(2\pi, ., \phi),\]

\noindent where $U^{(\epsilon)}$ is the solution to equation (\ref{bigU2}). Note that this is a bounded operator; 

\[ \|T\|_2 \leq \exp\{2\pi / \epsilon\}.\]

\noindent  Let $r(T) = \sup_{\lambda \in \sigma (T)} |\lambda|$, where $\sigma$ denotes spectrum.  By theorem VI.6 page 192 from Reed and Simon~\cite{RS}, $r(T) = \lim_{n \rightarrow +\infty} \|T^n\|^{1/n}_2$. It follows that $\exp\{2\pi / \epsilon\} \geq r(T) \geq \exp\{ -2\pi / \epsilon\}$. The operator $T$ is compact. Therefore, by the Riesz - Schauder theorem (Reed and Simon~\cite{RS} Theorem VI.15 page 203), $\sigma(T)$ is a discrete set having no limit points except perhaps for $0$. Therefore, there exists an eigenvalue $\lambda$ such that $|\lambda| = r(T)$ and this eigen value is of finite multiplicity. 

Let $\phi = \alpha + i\beta$ denote an eigenfunction, with $\|\phi\|_2 = \|\alpha\|_2 + \|\beta\|_2 = 1$, where $\alpha$ and $\beta$ are real functions, such that $\|T \phi\|_2 = r(T)$. Since $\phi$ (by hypothesis) {\em maximises} $\frac{\|T\phi\|_2}{\|  \phi\|_2}$, it is easy to see that  {\em both} $\alpha$ {\em and} $\beta$  maximise $\frac{\|T\phi\|_2}{\|\phi\|_2}$ and therefore that $\alpha$ is either non negative or non positive and $\beta$ is either non negative or non positive. The corresponding eigenvalue may be written as $\lambda = r(T)e^{i\theta}$, for $\theta \in [0,2\pi)$. Since both $\alpha$ and $\beta$ maximise the expression, it follows, using the notation $\langle \alpha, \beta \rangle = \frac{1}{2\pi}\int_0^{2\pi} \alpha(x)\beta(x)dx$, that 

\[ r^2\|\alpha\|_2^2 = \|T\alpha\|_2^2 = r^2\left(\cos^2(\theta) \|\alpha\|_2^2 + \sin^2(\theta) \|\beta\|_2^2 - 2\cos(\theta)\sin(\theta)\langle \alpha, \beta \rangle \right)
 \]

\noindent and

\[ r^2 \|\beta\|_2^2 = \|T\beta\|_2^2 = r^2 \left (\sin^2(\theta) \|\alpha\|_2^2 + \cos^2(\theta) \|\beta\|_2^2 + 2\sin(\theta)\cos(\theta) \langle \alpha, \beta \rangle \right ).\]

\noindent From this, $\theta = 0$ or $\pi$, from which it follows that $\lambda = r$; the eigenvalue is real and the eigenfunction may be taken as real and non negative.

\noindent Since $U^{(\epsilon)}(2\pi, ., \phi)$ is strictly positive for $\phi \in L^2$ non negative with $\|\phi\|_2 = 1$, it follows, since $\phi = \frac{1}{r}U^{(\epsilon)}(2\pi,.,\phi)$,  that $\phi$ is strictly positive. 

\noindent   Set

\[ u^{(\epsilon)}(0,.) = -\epsilon \frac{\partial}{\partial x}\log \phi.\]

\noindent This initial condition will provide periodic solutions to equation (\ref{youeps}). \vspace{5mm}

\noindent {\bf Part 2: Uniqueness for $\epsilon > 0$.} Suppose
that there are two periodic solutions,
$u^{(\epsilon,1)}$ and $u^{(\epsilon,2)}$, both of period $2\pi k$ in the time
variable. Set
$D = u^{(1)} - u^{(2)}$ and $S = u^{(1)} + u^{(2)}$. Firstly, by lemma
(\ref{ltwobound}), 

\[ \left (\frac{1}{2\pi}\int_0^{2\pi} S^2(t,x) dx \right )^{1/2}\leq
18\pi.\]

\noindent Next, for $\epsilon > 0$, $S(t,.) \in C^\infty ([0,2\pi])$ for all
$t \geq 0$. Since $S$ is periodic in both variables, this implies that 
there exists a constant $C < +\infty$ such that $\sup_{t,x} |S(t,x)| < C$.
 Let 
$(\alpha_{n}(t))_{n \in {\bf Z}}$ denote the Fourier coefficients of $D$; that is $D(t,x) =
\sum_{n}
\alpha_{n}(t)e^{inx}$. Let $\beta_{n} = \frac{\alpha_{n}}{in}$ for $n
\neq 0$ and $\beta_{0} \equiv 0$. Set $\tilde{D} = \sum_{n}
\beta_{n}(t)e^{inx}$. Then $\tilde{D}_x = D$. Note that $\tilde{D}$
satisfies

\[\tilde{D}_t = \frac{\epsilon}{2}\tilde{D}_{xx} - \frac{1}{2}S \tilde{D}_x.
\]

\noindent Set ${\cal L}(t,x) = \frac{\epsilon}{2}\frac{\partial^2}{\partial
x^2} - \frac{1}{2}S(t,x)\frac{\partial}{\partial x}$ and let $q(t;x,y)$ denote
the solution to 

\[ \left\{\begin{array}{l} \frac{\partial}{\partial t}
q(t;x,y) = {\cal L}(t,x)q(t;x,y) \\ q(0;x,y) = \delta_0(x-y).\end{array}\right.\]

\noindent Let ${\cal Q}$ denote the semigroup defined by

\[ {\cal Q} f(x) = \int q(2\pi;x,y)f(y) dy\]

\noindent and note that for all integer $N \geq 0$

\[ \tilde{D}(2\pi N, x) = {\cal Q}^N \tilde{D}(0,x).\]

 \noindent Because there exists a constant $C < +\infty$ such that
$\sup_{t,x}|S(t,x)| < +\infty$, it therefore follows by standard and straight forward results, for $\epsilon > 0$ that
$\inf_{0
\leq x \leq 2\pi}\inf_{0 \leq y \leq 2\pi} q(2\pi; x,y) > 0$, where the
inequality is strict. From this, it follows that 
${\cal Q}$ is the one step transition kernel of  an {\em ergodic} (discrete
time) time homogeneous Markov chain and hence a standard application of
the Ergodic theorem yields  that
$\lim_{N
\rightarrow +\infty} q(2N\pi, x,y) \rightarrow \tilde{q}(y)$, independent of
the initial condition $x$. Since $\tilde{D}$ is periodic, it follows that for all $N \geq 0$, 
$\tilde{D}(0,x) = \tilde{D}(2N\pi,x) = \int \tilde{q}(y)\tilde{D}(0,y)dy \equiv C$ where $C$ is a constant. Since $\int_0^{2\pi}
\tilde{D}(0,x) dx = 0$, it follows that $C \equiv 0$ and hence that
$\tilde{D} \equiv 0$.
 \vspace{5mm}

\noindent It follows that, for all $\epsilon > 0$, there exists a {\em
unique} initial condition for equation (\ref{youeps}) that provides solutions 
which are 
$2\pi$ periodic in the space variable and periodic in the time variable, and
that the solution is {\em unique}. In the time variable, the periodic solution
has period
$2\pi$. \vspace{5mm}

\noindent {\bf Uniqueness for $\epsilon = 0$} Consider equation
(\ref{youeps}) with $\epsilon = 0$. Suppose there exist two solutions
$2\pi$ periodic in the space variable and $2\pi k$ periodic in the time
variable. Denote the solutions by $u^{(1)}$ and $u^{(2)}$. Set $S = u^{(1)} +
u^{(2)}$. Then, as before, it is easy to see that there is a function
$\tilde{D}$ such that $\tilde{D}_x = u^{(1)} - u^{(2)}$ and such that
$\int_0^{2\pi}\tilde{D}(t,x)dx \equiv 0$, satisfying

\[   \tilde{D}_t = -\frac{1}{2}S \tilde{D}_x .\]

\noindent Let $Z$ denote the process defined by the relation

\[ Z_{s,t}(x) = x - \frac{1}{2}\int_s^t S(r,Z_{r,t}(x))dr.\]

\noindent Then

\[ \tilde{D}(t,x) =  \tilde{D}(0,Z_{0,t}(x)).\]

\noindent In particular, it holds for all non negative integer $N$ that 

\[ \tilde{D}(0,x)= \tilde{D}(2N\pi,x) =  \tilde{D}(0,Z_{0,2\pi N}(x)).\]

\noindent From lemma (\ref{krep}), it follows that $\tilde{D}(0,.)$ is constant on intervals $(x+(j-1)\frac{2\pi}{k},x+j\frac{2\pi}{k})$ for some positive integer $k$ and some $x \in [0,\frac{2\pi}{k})$. But since $u^{(1)} - u^{(2)}$ is bounded, it follows that $\tilde{D}(t,.)$ is Lipschitz in $x$ for each $t$ and hence constant (in the $x$ variable). Since $\int_0^{(2\pi)}\tilde{D}(t,x)dx \equiv 0$, it follows that $\tilde{D}(t,x) \equiv 0$ and hence uniqueness is established.

 It therefore follows
directly that, for $\epsilon = 0$, there is a unique initial condition for
equation (\ref{youeps}) which provides solutions to equation
(\ref{youeps}) that satisfy $\int_0^{2\pi} u(t,x) dx \equiv 0$, which are
$2\pi$ periodic in the space variable and periodic in the time variable. The
solution to equation (\ref{youeps}) with this initial condition is unique and
is $2\pi$ periodic in the time variable. Uniqueness of periodic solution to
equation (\ref{youeps}) with $\epsilon = 0$ satisfying $\int_0^{2\pi} u(t,x) dx \equiv 0$ follows.  
\vspace{5mm}

\noindent Using lemma (\ref{ltwobound}), there exists a weakly convergent
subsequence of
$u^{(\epsilon)}$ with limit $u$, which is
$2\pi$ periodic in both variables. Any such limit $u$ is a solution to
equation (\ref{youeps}) for $\epsilon = 0$. Since there
exists exactly one periodic solution to equation (\ref{youeps}) for
$\epsilon = 0$, it follows that the sequence $u^{(\epsilon)}$ converges in
$L^2$, in the relative weak topology, to
$u$.  The proof of lemma (\ref{periodexiuni}) is now complete.
\qed\vspace{5mm} 

\begin{Lmm} \label{endgame} Let $u$ denote the periodic solution to equation
(\ref{youeps}) with $\epsilon = 0$. Let $\tilde{\theta}$ solve

\begin{equation}\label{theta} \left\{\begin{array}{l} \dot{\tilde{\theta}} =
-u(-t,
\tilde{\theta})\\
\tilde{\theta} (0,x) = x
\end{array}\right. 
\end{equation} 

\noindent Then, for all $x \in [0,2\pi)$, either 

\[ \limsup_{t \rightarrow +\infty}  |(\tilde{\theta} (t,x) + t - \sin(t))
\mod(2\pi) | = 0,\]

\noindent or

\[ \limsup_{t \rightarrow +\infty} |(\tilde{\theta} (t,x) - \pi - t -
\sin(t))\mod(2\pi) | = 0.\]

\end{Lmm}\vspace{5mm}

\noindent {\bf Proof of lemma (\ref{endgame})} Firstly, lemma (\ref{krep}) gives that there exists a positive integer $k$ such that for all $x \in [0,2\pi)$ and $y \in [0,2\pi)$, 

\begin{equation}\label{endnice}\lim_{t \rightarrow +\infty}|k (\tilde{\theta}(t,x) - \tilde{\theta}(t,y))| = 0.
\end{equation}

\noindent Note that $\tilde{\theta}(t,x) = \theta(-t,x)$, where $\theta$ solves 

\[ \left\{\begin{array}{l}\dot{\theta}(t,x) = u(t,\theta(t,x)),\\ \theta(0,x) = x. \end{array}\right. \]

\noindent Set 

\begin{equation}\label{essdef} {\cal S} := \{(q,p) | q \in {\bf S}^1, p =
-u(0,q)\}\end{equation}

\noindent and set $T: {\cal S} \rightarrow {\cal S}$, defined such that

\begin{equation}\label{teedef} T(q,p) = (\tilde{\theta}(2\pi, q), -u(2\pi,
\tilde{\theta} (2\pi,q)).\end{equation}

\noindent  Then it is clear by lemma (\ref{disappear}), from the construction of the solution to the inviscid Burgers equation, given by the method of characteristics, described in  section (\ref{downward}) that for all $n \geq 0$,
$T^{(n+1)} {\cal S}
\subseteq T^{(n)}{\cal S}$. Set 

\begin{equation}\label{tildeess}\tilde{\cal S} = \cap_{n \geq 1} \overline{T^{(n)}{\cal S}}.\end{equation}

\noindent Then
it is clear that
$\tilde{\cal S}$ is non empty and that $T\tilde{\cal S} = \tilde{\cal S}$. Furthermore, from equation (\ref{endnice}), it follows that there exists a positive integer $k$ such that 

\[\tilde{\cal S} = \{(q + \frac{2\pi j}{k}, -u(0,q + \frac{2\pi j}{k})|j = 0,1, \ldots, k\}.\]

\noindent It is now shown that $\tilde{\cal S}$ consists of a single point; either
$\tilde{\cal S} = \{(0,0)\}$, or $\tilde{\cal S} = \{(\pi,2 ) \}$. \vspace{5mm}

\noindent By construction, the points $q + \frac{2\pi j}{k}$ provide initial conditions that give periodic solutions to the equation

\[ \dot{\theta} = u(t,\theta).\]

\noindent Note that $\theta (-t,.) = \tilde{\theta}(t,.)$.   From the construction, the trajectories $\theta$ with these initial conditions survive for all time in the construction of the inviscid Burgers equation described in section (\ref{downward}). They are not absorbed into a downward jump and therefore they satisfy

\[ \left\{ \begin{array}{l} \ddot{\theta} = \sin (\theta + \sin(s))\\ \theta(0) = q + \frac{2\pi j}{k}, \quad \dot{\theta}(0) = u(0,q + \frac{2\pi j}{k}), \quad j=0,1,\ldots, k-1.
\end{array}
 \right.\]

\noindent Since these trajectories do not intersect, it follows from theorem (\ref{ONEsol}) that $k=1$ and that EITHER $\tilde{S} = (0,0)$ OR $\tilde{S} = (\pi,2)$. From this, it follows that EITHER

\[  \limsup_{t \rightarrow +\infty}  |(\tilde{\theta}(t,x) + t -
\sin(t))
\mod(2\pi) | = 0 \qquad \forall x \in {\bf S}^1,\]

\noindent OR

\[  \limsup_{t \rightarrow +\infty} |(\tilde{\theta} (t,x) - \pi - t -
\sin(t))\mod(2\pi) | = 0 \qquad \forall x \in {\bf S}^1\]

\noindent and  lemma (\ref{endgame})is proved. 
 \qed \vspace{5mm}

\noindent Let $\phi^{(\epsilon)}(.)$ denote the function such that  
$u^{(\epsilon)}(0,.) :=
\phi_x^{(\epsilon)}(.)$ gives a periodic solution to equation (\ref{youeps}).
Then $\phi^{(\epsilon)}_x$ converges in $L^2$, in the relative weak topology
to   $\phi_x$ and, furthermore, $\sup_{0 \leq \epsilon \leq 1} \|\phi_x^{(\epsilon)} \|_2 \leq 6\pi + \sqrt{6(1+\pi)}$ by lemma (\ref{ltwobound}). It
follows by the Ascoli Arzela lemma that  
$\phi^{(\epsilon)}$ has a limit
$\phi$ such that $\lim_{\epsilon
\rightarrow 0}\sup_{0 \leq x \leq 2\pi} |\phi^{(\epsilon)} (x) -
\phi (x) | = 0$.
  \vspace{5mm}

\noindent {\bf Completing the Counter Example} The arguments given in the 
previous sections are all standard and yield the following: let
$u^{(\epsilon)}$ denote the unique periodic solution to equation (\ref{youeps}),
then the sequence
$u^{(\epsilon)}$ is  convergent in $L^2$, in the relative weak topology, to
a limit $u$, which provides the unique periodic solution to equation (\ref{youeps})
with $\epsilon = 0$. Firstly,  

\[ u (t,x) = -\lim_{\epsilon \rightarrow 0}\frac{\partial}{\partial x} 
\log E_{{\bf P}}
\left[ e^{-\frac{1}{\epsilon }\left(\phi^{(\epsilon )}(x +
\sqrt{\epsilon } w_t) + \int_0^t \cos(\sin(t-s))ds - \int_0^t
\cos(\sin(t-s) + x + \sqrt{\epsilon }w(s))ds \right )}\right],\]

\noindent yielding

\[u(t,x) = \frac{\partial}{\partial x}\inf_{\xi : \xi(t) = x} \left \{ \phi(\xi(t)) +
\int_0^t \cos (\sin s) ds - \int_0^t \cos (\xi(s) + \sin s) ds \right \},\]

\noindent by Varadhan's theorem (\ref{varad}) and the results of section (\ref{varadhansect}). Next, the function $u$ has representation 

\begin{equation}\label{firstrep} u(t,x) = \dot{\eta}(t, \eta^{-1}(t,x)),
\end{equation}

\noindent where $\eta$ solves 

\begin{equation}\label{rep2} \left\{\begin{array}{l} \ddot{\eta} = \sin (\eta
+
\sin t)
\\
\eta(0,x) = x,
\;
\dot{\eta}(0,x) = \phi_x(x).\end{array}\right.\end{equation}

\noindent This may be rewritten as

\begin{equation}\label{firstrep2} u(t,x) = \dot{\xi}^{(t,x)}(t)
\end{equation}

\noindent where $\xi$ solves

\begin{equation}\label{rep22}\left\{\begin{array}{ll}\ddot{\xi}^{(t,x)}(s) = \sin
(\xi^{(t,x)} (s)+ \sin(s)) & 0 \leq s \leq t \\ \xi^{(t,x)}(t) = x, \quad \dot{\xi}^{(t,x)}(0) = \phi_x
(\xi^{(t,x)}(0)).\end{array}\right. \end{equation}

\noindent  Furthermore, it was also
established that the solutions
$\xi$ used in representation (\ref{firstrep2}), which solve equation
(\ref{rep22}), {\em  minimise} the action functional 

\[ {\cal A}(\xi;t,x) := 
\left
\{
\phi(\xi(t)) +
\int_0^t \cos (\sin s) ds - \int_0^t \cos (\xi(s) + \sin s) ds \right \}\]

\noindent 
subject to the constraint that $\xi(t) = x$.

\noindent  Let $\tilde{\theta}(s,x) = \xi^{(0,x)}(-s)$, then 

\[ \left\{\begin{array}{l} \ddot{\tilde{\theta}} = \sin(\tilde{\theta} - \sin s)\\ \tilde{\theta}(0,x)
= x,
\; \dot{\tilde{\theta}}(0,x) = -u(0,x).\end{array}\right. \]

\noindent Note that 

\[ \dot{\tilde{\theta}} = -u(-t,\tilde{\theta})\]

\noindent for all $t \geq 0$.  Recall the definition of ${\cal S}$ given in
equation (\ref{essdef})  and
recall the definition of $T$ given in equation (\ref{teedef}). Recall the definition of $\tilde{\cal S}$ given by equation (\ref{tildeess}) and recall that $\tilde{\cal S} = \{(0,0)\}$ or $\{(\pi, 2)\}$.   The points of $\tilde{\cal S}$ give initial conditions $(\tilde{\theta}(0), \dot{\tilde{\theta}}(0))$ that yield 
{\em non intersecting} periodic solutions to the equation

\[ \ddot{\tilde{\theta}} = \sin (\tilde{\theta} - \sin (t)).\]

\noindent Lemma (\ref{ONEsol}) showed that, modulo $2\pi$, there
existed only two such trajectories; 

\[ \tilde{\theta} (t) = -t + \sin (t) \]

\noindent and

\[ \tilde{\theta}(t) = \pi + t + \sin(t).\]

\noindent Since these intersect, therefore $\tilde{\cal S}$ can have at most one
point.  The point
$(q,p)$, where $(q,-p) \in \tilde{\cal S}$ represents the initial condition
$(\eta(0),\dot{\eta}(0)) = (q,p)$ for trajectories solving equation
(\ref{rep2}) that survive for all time in the construction (\ref{firstrep})
of   periodic solutions to the inviscid Burgers equation. It follows
that $\tilde{\cal D}$ contains exactly one point,   either $(0,0)$ or
$(\pi, 2)$.  It follows that there is exactly one Euler Lagrange
trajectory in the construction of the solutions to the inviscid
Burgers' equation that survives for all time and that  the trajectory is
EITHER 

\[ \xi(t) = t - \sin(t) \]

\noindent OR  

\[ \xi(t) = \pi - 2t - \sin (t).\] 

\noindent But lemma (\ref{EZLMM}) proves that neither of these minimise
the associated action functional for $t \geq 2\pi$. It shows that for any periodic solution to
the inviscid Burgers' equation under consideration, there {\em necessarily} exist trajectories
of the associated Euler Lagrange equation used in the construction which do
not minimise the action functional. It follows that the minimising
trajectories do not yield periodic solutions to the inviscid Burgers'
equation.
\vspace{5mm} 

\noindent For sufficiently large $t$, the periodic viscosity solution to the equation

\[ u_t + \frac{1}{2}(u^2)_x = \sin (x + \sin t)\]

\noindent therefore may {\em not} be constructed using  the associated
Euler Lagrange equations which {\em minimise} the action functional. 

A counterexample has therefore been given  to the assertion
that the viscosity solution to the inviscid Burgers' equation in the case of
smooth, bounded space / time periodic potentials may always be
constructed from the Euler Lagrange trajectories that minimise the
associated action functional. 

On the other hand, a full proof of this assertion 
exists, and has been outlined in this note, if Tychonov's theorem is
assumed. The negation of Tychonov's theorem implies the negation of the
Axiom of Choice. It follows that Tychonov's theorem, and hence the Axiom
of Choice lead to contradictions and are therefore inadmissible in
mathematical analysis.   
\qed \vspace{5mm}

\noindent Electronic mail address for correspondence: {\tt jonob@mai.liu.se}\vspace{5mm}
  
\setlength{\baselineskip}{2ex}

\end{document}